%% file: hunger-game.tex
\documentclass[11pt]{amsart}
\usepackage{preamble}
\usepackage{hyperref}

\begin{document}

\title[A greedy chip-firing game]{A greedy chip-firing game}
\author[Rupert Li]{Rupert Li}
\address{Massachusetts Institute of Technology}
\email{\href{mailto:rupertli@mit.edu}{{\tt rupertli@mit.edu}}}
\author[James Propp]{James Propp}
\address{University of Massachusetts Lowell, Department of Mathematical Sciences}
\email{\href{mailto:jamespropp@gmail.com}{{\tt jamespropp@gmail.com}}}
\date{\today}

\begin{abstract}
We introduce a deterministic analogue of Markov chains 
that we call the hunger game.  
Like rotor-routing, the hunger game
deterministically mimics the behavior of both recurrent Markov chains 
and absorbing Markov chains.
In the case of recurrent Markov chains with finitely many states, 
hunger game simulation concentrates around the stationary distribution
with discrepancy falling off like $N^{-1}$,
where $N$ is the number of simulation steps;
in the case of absorbing Markov chains with finitely many states, 
hunger game simulation also exhibits concentration 
for hitting measures and expected hitting times
with discrepancy falling off like $N^{-1}$
rather than $N^{-1/2}$.
When transition probabilities in a finite Markov chain are rational,
the game is eventually periodic;
the period seems to be the same for all initial configurations
and the basin of attraction appears to tile the configuration space 
(the set of hunger vectors) by translation,
but we have not proved this.
\end{abstract}

\maketitle

\import{}{introduction}
\import{}{preliminaries}
\import{}{boundedness}
\import{}{termination}

\import{}{stationary}
\import{}{hitting}
\import{}{recurrence}
\import{}{conclusion}

\section*{Acknowledgments}
The authors thank Grant Barkley, Darij Grinberg,
Sam Hopkins, Lionel Levine, Alex Postnikov, and Tom Roby
for various helpful suggestions made during the course of this research.
We are grateful to the anonymous referees who made many suggestions that improved this article.

\bibliographystyle{plain}
\bibliography{ref}

\end{document}

%% file: introduction.tex
\section{Introduction}\label{section: introduction}
In the 1970s, Arthur Engel \cite{engel1975probabilistic,engel1976does} 
introduced a ``stochastic abacus''
that deterministically mimics many aspects of the behavior
of finite-state Markov chains with rational probabilities.
Unaware of Engel's work, various mathematicians and physicists
invented and studied the abelian sandpile model \cite{dhar1990self} 
and the chip-firing game \cite{bjorner1991chip}
which in many respects embody the same core idea as Engel's abacus
but with different motivations.
For books surveying chip-firing we refer readers to 
\cite{klivans2018mathematics,corry2018divisors}.

In the 2010s, inspired by Engel's work, the second author of this paper 
in collaboration with Ander Holroyd 
\cite{holroyd2010rotor} 
introduced a different way to deterministically mimic Markov chains
via ``rotor-routing'' (unaware that physicists were already studying
the process under the name ``the Eulerian walkers model'' 
\cite{priezzhev1996eulerian}).
The emphasis of much of this work on the rotor-router model,
inspired by discrepancy theory and quasi-Monte Carlo methods,
was on the fidelity of the deterministic process to the associated Markov chain.
Specifically, it was shown that for certain 
asymptotically-defined quantities associated with Markov chains
(e.g., the proportion of the time that the chain spends in a specific state),
rotor-router simulation has the same limiting behavior as the Markov chain,
typically with faster convergence than 
the Markov chain itself would typically exhibit.
For a comprehensive background on chip-firing and rotor-routing and 
the relationship between them, we refer readers to \cite{holroyd2008chip}.

Here we discuss another way to derandomize Markov chains
which we call the {\em hunger game}.
It can be applied to any discrete-state discrete-time Markov chain,
whether or not the transition probabilities are rational.
(Rotor-routing can be extended to this regime --- 
see the discussion of stack walks in \cite{holroyd2010rotor} ---
though we know of no way to extend chip-firing in this direction.)
A key difference between rotor-routing and the hunger game
is that the frequency with which state $i$ is followed by $j$
in rotor-router simulation converges to the transition probability $P_{ij}$;
this is not the case for the hunger game.
That is, the hunger game does not exhibit fidelity 
with regard to transition-frequencies. 
However, we believe this may be a virtue rather than a vice,
as we will explain in the concluding section.

In \cref{section: preliminaries} we define the hunger game
and the chip addition operators associated with it.
In \cref{section: boundedness hunger}
we prove fundamental results on the boundedness of hunger,
and in \cref{section: termination} we prove basic results 
on the behavior of chip addition operators.
In \cref{section: stationary distributions} we prove 
the main result of this paper, demonstrating that 
the normalized firing vector of a hunger game process converges 
to the unique stationary distribution of an irreducible Markov 
chain with a discrepancy bound inversely proportional 
to the number of time steps.
We apply this result to show how the hunger game process 
can calculate hitting probability distributions, escape probabilities, 
expected absorption times, and expected return times 
in \cref{section: hitting probabilities and hitting times}\,.
In \cref{section: recurrence} we focus on finite Markov chains
with rational transition probabilities.
We introduce the notion of a recurrent hunger vector 
(a vector that returns to itself under the hunger game) 
and study the properties of the basin of attraction
(the set of recurrent vectors).
In the case where all transition probabilities are rational,
we prove the zero vector is always recurrent and determine its period, 
and conjecture that all the periods for a given Markov chain are the same.
We conclude with \cref{section: conclusion}\,,
comparing the discrepancies of 
the rotor-router model and the hunger game.

%% file: preliminaries.tex
\section{Preliminaries}\label{section: preliminaries}
See \cite{kemenysnell1983} or \cite{norris1998markov}
for basic facts about Markov chains.
Throughout this paper, except where otherwise noted, 
Markov chains are assumed to have a finite state space
indexed by positive integers $1,\dots,n$ for some $n$; 
when we consider Markov chains with countably infinite state spaces, 
we will assume that for each state $i$
there are only finitely many states $j$ such that 
the transition probability $P_{ij}$ from state $i$ to state $j$ is positive.
We can represent a Markov chain by a weighted directed graph $G=(V,E)$ 
whose vertices $v_i \in V$ are the allowed states 
and the weight of the directed edge $(v_i,v_j) \in E$ is $P_{ij}$.
For brevity, we will use the words state and vertex to refer to 
both the state in the Markov chain and the corresponding vertex in $G$,
passing back and forth between the abstract Markov chain
and its concrete embodiment as a random walk on $G$.
We say that state $i$ is an \textbf{absorbing state} if $P_{ii}=1$
(equivalently if $G$ has a loop at $v_i$ with weight 1).
Notice that for every vertex $v \in G$, 
the sum of the weights of all edges $(v,w)$ is 1. 

Let $H=P-I$, where $I$ is the $n$-by-$n$ identity matrix.
Let $H_i$, $P_i$, and $I_i$ denote the $i$th rows 
of the matrices $H$, $P$, and $I$ respectively.
Observe that $I_i = e_i$ (the $i$th unit vector)
and that $-H$ is the Laplacian of the graph.
The Markov chain admits at least one \textbf{stationary measure} $\ppi$
for which the vector $v = [\ppi(1),\dots,\ppi(n)]$ satisfies $vP=v$
(this follows from the assumption that the Markov chain is finite;
see, e.g., \cite{norris1998markov}),
so 1 is an eigenvalue of $P$ and 0 is an eigenvalue of $H$.
All eigenvalues of $P$ have magnitude at most 1.
When the Markov chain is \textbf{irreducible}
(that is, when every state can be reached from every other state
in some finite number of steps), 1 is a simple eigenvalue of $P$
and 0 is a simple eigenvalue of $H$,
so that $H$ has rank $n-1$,
and the rows of $H$, taken with integer coefficients,
generate an $n-1$-dimensional sublattice of the 
space of vectors with entries summing to zero;
in this case, there is a unique stationary measure $\ppi$
satisfying $\ppi(1)+\cdots+\ppi(n) = 1$.

We now informally introduce the hunger game on the weighted directed graph $G$
by comparing it to the chip-firing model and the rotor-router model
before offering a more technical definition.

The ``goodness'' of a deterministic analogue of a random process
can be assessed by the notion of discrepancy.
If some numerical characteristic of the deterministic process
converges to a corresponding numerical characteristic of the random process
as simulation time goes to infinity,
one can try to determine the rate of convergence.

The simplest deterministic analogues of Markov chains
were invented by Engel \cite{engel1975probabilistic,engel1976does}, 
under the name \textbf{the stochastic abacus}
(though the term \textbf{chip-firing} is more common nowadays).
Suppose that all the transition probabilities $P_{ij}$ are rational,
and that we have positive integers $d_1,\dots,d_n$
such that $d_i P_{ij}$ is an integer for all $i,j$.
Assume that the Markov chain is irreducible.
Define a \textbf{chip-configuration} as 
an $n$-tuple $\mathbf{c} = (\mathbf{c}_1,\dots,\mathbf{c}_n)$;
say that a chip-configuration is \textbf{stable}
if $\mathbf{c}_i < d_i$ whenever $i$ is a non-absorbing state.
If $\mathbf{c}_i \geq d_i$, then we obtain another chip-configuration 
by \textbf{firing at $i$}, replacing $\mathbf{c}$ by
$\mathbf{c'} = \mathbf{c} + d_i H_i$;
if we represent the Markov chain by drawing a graph with vertices $v_1,\dots,v_n$
and we represent the chip-configuration $\mathbf{c}$ 
by putting $\mathbf{c}_j$ chips at $v_j$ for all $j$,
then firing means sending $d_i P_{ij}$ chips 
from $i$ to $j$ for each $j \neq i$.
Chip-firing can be used to find the stationary probability measure 
of an irreducible Markov chain as follows.
Put sufficiently many chips on the state-graph so that 
stabilization is impossible no matter how many firings are performed, 
and start performing firings however one wishes;
when we encounter a chip-configuration we have seen before
(as must happen eventually), the vector that records the number of times
each state has fired will be a stationary vector.

\begin{example}\label{example: 3 chip}
Suppose we have the Markov chain given by the Markov matrix
\begin{align*}
    \begin{bmatrix}
        \frac{1}{2} & \frac{1}{2} & 0 \\
        \frac{1}{2} & 0 & \frac{1}{2} \\
        0 & \frac{1}{2} & \frac{1}{2}
    \end{bmatrix},
\end{align*}
representing a doubly-reflecting random walk.
Its corresponding graph is shown in \cref{fig:ex2G}\,.
\begin{figure}[htbp]
    \centering
    \begin{tikzpicture}[scale=0.6,font=\normalsize,baseline,thick]
        \foreach \x in {1,...,3} {
            \filldraw[color=black,fill=green!30,thick] (4*\x,0) circle (16pt);
            \node at (4*\x,0) {$v_{\x}$};
        }
        \foreach \x in {2,3} {
            \draw [->,>=latex] plot [smooth,tension=1] coordinates {(4*\x-0.866*0.7,0.4*0.7) (4*\x-2,0.7) (4*\x-4+0.866*0.7,0.4*0.7)};
        }
        \foreach \x in {1,2} {
            \draw [->,>=latex] plot [smooth,tension=1] coordinates {(4*\x+0.866*0.7,-0.4*0.7) (4*\x+2,-0.7) (4*\x+4-0.866*0.7,-0.4*0.7)};
        }
        \node at (6,1) {0.5};
        \node at (10,1) {0.5};
        \node at (6,-1) {0.5};
        \node at (10,-1) {0.5};
        \draw [->,>=latex] plot [smooth,tension=5] coordinates {(4-0.866*0.7,0.5*0.7) (4-1.5,0) (4-0.866*0.7,-0.5*0.7)};
        \node at (2,0) {0.5};
        \draw [->,>=latex] plot [smooth,tension=5] coordinates {(12+0.866*0.7,-0.5*0.7) (12+1.5,0) (12+0.866*0.7,0.5*0.7)};
        \node at (14,0) {0.5};
    \end{tikzpicture}
    \caption{A graph $G$ corresponding to a doubly-reflecting random walk.}
    \label{fig:ex2G}
\end{figure}
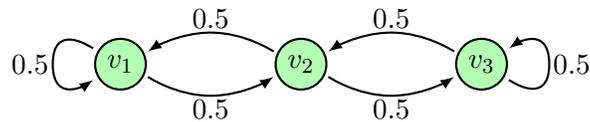
If we place 1 chip at state 1, 2 chips at state 2, and 1 chip at state 3,
then state 2 is unstable, so we may fire at 2,
turning the chip-configuration (1,2,1) into the chip-configuration (2,0,2).
Then firing at 1 and at 3 brings us back to (1,2,1).
The vector that records the number of times each state fired
is (1,1,1), which is indeed a stationary vector for this Markov chain.
\end{example}

The rotor-router model~\cite{holroyd2010rotor}
is a different scheme for imitating Markov chains deterministically.
Assume as above that the Markov chain
has rational transition probabilities, with $d_i$ as above.
Represent the Markov chain using a directed graph
with $d_i P_{ij}$ parallel edges from $v_i$ to $v_j$,
so that $v_i$ has outdegree $d_i$.
Each vertex distributes arriving chips along its outgoing edges in a cyclic manner, 
not sending a chip along any edge for a second time
until it has sent a chip along every edge at least once,
and thereafter always sending the next chip along 
the edge along which it has sent a chip least recently.
Inasmuch as the vertex with the chip gets to decide where the chip goes next,
we call this ``supply-side'' management of the chip's movement.
Assume that the Markov chain is irreducible.
It can be shown that once the chip enters an infinite loop
(as must happen eventually),
the fraction of the time that the chip spends at vertex $v_i$
is proportional to the steady-state $\ppi(i)$.

\begin{example}\label{example: 3 rotor}
We use the same Markov chain as~\cref{example: 3 chip}.
Suppose we start with the chip at $v_1$
and begin the game by sending the chip to $v_2$, then $v_3$, then $v_3$ again.
Since at this point the chip has already traveled 
along the edge sending $v_3$ to $v_3$,
it now travels along the edge from $v_3$ to $v_2$.
As the chip has already gone from $v_2$ to $v_3$,
the rotor-router protocol dictates
that it must now travel along the other edge from $v_2$ and go to $v_1$.
Similarly, as it has already gone from $v_1$ to $v_2$,
now it must travel from $v_1$ to $v_1$.
Thereafter the process cycles forever.
Since within each cycle
the chip spends 2 steps at each vertex,
(2,2,2) is a stationary vector.
\end{example}

The \textbf{hunger game} introduced in this article
can be seen as a ``demand-side" management system: 
each vertex has a ``hunger" for chips, 
determined by its expectation of receiving chips 
from neighboring vertices that have been previously visited.
When a vertex $v$ receives a chip, the neighboring vertices' hunger 
increases in accordance to the transition probabilities from $v$ 
to those vertices, and $v$ sends its chip to the vertex with highest hunger
(regardless of whether that chip is a neighbor of $v$).

Now we give a more formal definition of the \textbf{hunger game}.
It is simplest to start with the situation
in which the chain runs forever without restarts
(in contrast to chains that will be restarted 
when they enter an absorbing state).
We also start with the case in which the state space is finite,
with $|V| = n$, deferring discussion of infinite-state spaces until later.
The \textbf{hunger vector} $\mathbf{h}\in\R^{n}$ represents 
the hunger at each vertex in $V$;
we will sometimes call it the \textbf{hunger state} to emphasize 
its interpretation as a state of the hunger game system.
At each step, whichever vertex has the highest hunger receives the chip; 
if vertex $v_i$ receives the chip, then the hunger vector $\mathbf{h}$ 
is updated by adding $H_i = P_i - I_i$ to it, corresponding to 
the increase in vertices' hunger from the presence of this chip at $v_i$ 
but also the satiation of $v_i$ after receiving a chip.
If multiple states are tied for the highest hunger,
we break the tie by choosing the lowest-indexed such vertex.
Since each row of $H$ has entries summing to 0,
total hunger never changes.

Each relocation of the chip is referred to as \textbf{firing} the chip;
specifically, when the chip is relocated at $i$,
we say the chip fires to $i$.
Unlike rotor-router or chip-firing, 
under the hunger game rules a chip does not necessarily have to be fired 
to a vertex adjacent to its current location; see \cref{example: 3 hunger}\,.
In fact, the determination of the chip's next location depends 
only on the hunger state $\mathbf{h}$ and not the current location of the chip.
The step can be described purely in terms of
the matrix $H$ and the vector $\mathbf{h}$
without any reference to chips, via the rule
``Replace $\mathbf{h}$ by $\mathbf{h}' = \mathbf{h} + H_i$
where $i$ maximizes $\mathbf{h}_i$,
choosing the smallest such $i$ in the event of a tie.'' 

\begin{example}\label{example: 3 hunger}
We use the same Markov chain as~\cref{example: 3 chip}.
Starting with $\mathbf{h}=0$, as shown in \cref{subfig:ex2init}\,, 
regardless of the initial location of the chip, we fire the chip to $v_1$ 
under the tie-breaking rule, 
yielding $\mathbf{h}=\left[-\frac{1}{2},\frac{1}{2},0\right]$ 
as shown in \cref{subfig:ex2fire1}\,.
After this, we fire to $v_2$, as shown in \cref{subfig:ex2fire2}\,, 
and then fire to $v_3$, as shown in \cref{subfig:ex2fire3}\,.
Notice that we have returned back to 
the initial hunger state $\mathbf{h}=\mathbf{0}$, 
so this process repeats, 
visiting states 1, 2, 3, then back to 1, and so on.
Notice that we fire the chip from 3 to 1
even though $P_{31}=0$.
\begin{figure}[htbp]
    \centering
    \begin{subfigure}{\textwidth}
        \centering
        \begin{tikzpicture}[scale=0.6,font=\normalsize,baseline,thick]
            \foreach \x/\xcolor in {1/white,2/white,3/white} {
                \filldraw[color=black,fill=\xcolor,thick] (4*\x,0) circle (16pt);
            }
            \foreach \x/\xtext in {1/0,2/0,3/0} {
                \node at (4*\x,0) {\xtext};
            }
            \foreach \x in {2,3} {
                \draw [->,>=latex] plot [smooth,tension=1] coordinates {(4*\x-0.866*0.7,0.4*0.7) (4*\x-2,0.7) (4*\x-4+0.866*0.7,0.4*0.7)};
            }
            \foreach \x in {1,2} {
                \draw [->,>=latex] plot [smooth,tension=1] coordinates {(4*\x+0.866*0.7,-0.4*0.7) (4*\x+2,-0.7) (4*\x+4-0.866*0.7,-0.4*0.7)};
            }
            \node at (6,1) {0.5};
            \node at (10,1) {0.5};
            \node at (6,-1) {0.5};
            \node at (10,-1) {0.5};
            \draw [->,>=latex] plot [smooth,tension=5] coordinates {(4-0.866*0.7,0.5*0.7) (4-1.5,0) (4-0.866*0.7,-0.5*0.7)};
            \node at (2,0) {0.5};
            \draw [->,>=latex] plot [smooth,tension=5] coordinates {(12+0.866*0.7,-0.5*0.7) (12+1.5,0) (12+0.866*0.7,0.5*0.7)};
            \node at (14,0) {0.5};
        \end{tikzpicture}
        \caption{The initial hunger state $\mathbf{h}=\mathbf{0}$.}
        \label{subfig:ex2init}
    \end{subfigure}
    \begin{subfigure}{\textwidth}
        \centering
        \begin{tikzpicture}[scale=0.6,font=\normalsize,baseline,thick]
            \foreach \x/\xcolor in {1/blue!30,2/yellow!30,3/white} {
                \filldraw[color=black,fill=\xcolor,thick] (4*\x,0) circle (16pt);
            }
            \foreach \x/\xtext in {1/-0.5,2/0.5,3/0} {
                \node at (4*\x,0) {\xtext};
            }
            \foreach \x in {2,3} {
                \draw [->,>=latex] plot [smooth,tension=1] coordinates {(4*\x-0.866*0.7,0.4*0.7) (4*\x-2,0.7) (4*\x-4+0.866*0.7,0.4*0.7)};
            }
            \foreach \x in {1,2} {
                \draw [->,>=latex] plot [smooth,tension=1] coordinates {(4*\x+0.866*0.7,-0.4*0.7) (4*\x+2,-0.7) (4*\x+4-0.866*0.7,-0.4*0.7)};
            }
            \node at (6,1) {0.5};
            \node at (10,1) {0.5};
            \node at (6,-1) {0.5};
            \node at (10,-1) {0.5};
            \draw [->,>=latex] plot [smooth,tension=5] coordinates {(4-0.866*0.7,0.5*0.7) (4-1.5,0) (4-0.866*0.7,-0.5*0.7)};
            \node at (2,0) {0.5};
            \draw [->,>=latex] plot [smooth,tension=5] coordinates {(12+0.866*0.7,-0.5*0.7) (12+1.5,0) (12+0.866*0.7,0.5*0.7)};
            \node at (14,0) {0.5};
        \end{tikzpicture}
        \caption{$\mathbf{h}$ after firing to $v_1$, shown in blue. States with updated hungers are shown in yellow.}
        \label{subfig:ex2fire1}
    \end{subfigure}
    \begin{subfigure}{\textwidth}
        \centering
        \begin{tikzpicture}[scale=0.6,font=\normalsize,baseline,thick]
            \foreach \x/\xcolor in {1/yellow!30,2/blue!30,3/yellow!30} {
                \filldraw[color=black,fill=\xcolor,thick] (4*\x,0) circle (16pt);
            }
            \foreach \x/\xtext in {1/0,2/-0.5,3/0.5} {
                \node at (4*\x,0) {\xtext};
            }
            \foreach \x in {2,3} {
                \draw [->,>=latex] plot [smooth,tension=1] coordinates {(4*\x-0.866*0.7,0.4*0.7) (4*\x-2,0.7) (4*\x-4+0.866*0.7,0.4*0.7)};
            }
            \foreach \x in {1,2} {
                \draw [->,>=latex] plot [smooth,tension=1] coordinates {(4*\x+0.866*0.7,-0.4*0.7) (4*\x+2,-0.7) (4*\x+4-0.866*0.7,-0.4*0.7)};
            }
            \node at (6,1) {0.5};
            \node at (10,1) {0.5};
            \node at (6,-1) {0.5};
            \node at (10,-1) {0.5};
            \draw [->,>=latex] plot [smooth,tension=5] coordinates {(4-0.866*0.7,0.5*0.7) (4-1.5,0) (4-0.866*0.7,-0.5*0.7)};
            \node at (2,0) {0.5};
            \draw [->,>=latex] plot [smooth,tension=5] coordinates {(12+0.866*0.7,-0.5*0.7) (12+1.5,0) (12+0.866*0.7,0.5*0.7)};
            \node at (14,0) {0.5};
        \end{tikzpicture}
        \caption{$\mathbf{h}$ after firing to $v_2$, shown in blue. States with updated hungers are shown in yellow.}
        \label{subfig:ex2fire2}
    \end{subfigure}
    \begin{subfigure}{\textwidth}
        \centering
        \begin{tikzpicture}[scale=0.6,font=\normalsize,baseline,thick]
            \foreach \x/\xcolor in {1/white,2/yellow!30,3/blue!30} {
                \filldraw[color=black,fill=\xcolor,thick] (4*\x,0) circle (16pt);
            }
            \foreach \x/\xtext in {1/0,2/0,3/0} {
                \node at (4*\x,0) {\xtext};
            }
            \foreach \x in {2,3} {
                \draw [->,>=latex] plot [smooth,tension=1] coordinates {(4*\x-0.866*0.7,0.4*0.7) (4*\x-2,0.7) (4*\x-4+0.866*0.7,0.4*0.7)};
            }
            \foreach \x in {1,2} {
                \draw [->,>=latex] plot [smooth,tension=1] coordinates {(4*\x+0.866*0.7,-0.4*0.7) (4*\x+2,-0.7) (4*\x+4-0.866*0.7,-0.4*0.7)};
            }
            \node at (6,1) {0.5};
            \node at (10,1) {0.5};
            \node at (6,-1) {0.5};
            \node at (10,-1) {0.5};
            \draw [->,>=latex] plot [smooth,tension=5] coordinates {(4-0.866*0.7,0.5*0.7) (4-1.5,0) (4-0.866*0.7,-0.5*0.7)};
            \node at (2,0) {0.5};
            \draw [->,>=latex] plot [smooth,tension=5] coordinates {(12+0.866*0.7,-0.5*0.7) (12+1.5,0) (12+0.866*0.7,0.5*0.7)};
            \node at (14,0) {0.5};
        \end{tikzpicture}
        \caption{$\mathbf{h}$ after firing to $v_3$, shown in blue. States with updated hungers are shown in yellow.}
        \label{subfig:ex2fire3}
    \end{subfigure}
    \caption{The hunger game on a doubly-reflecting random walk.}
    \label{fig:ex2H}
\end{figure}
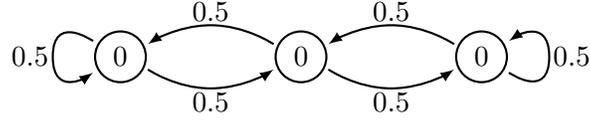
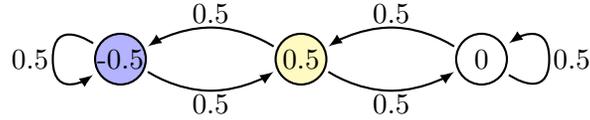
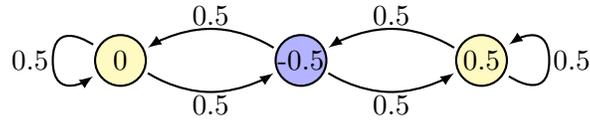
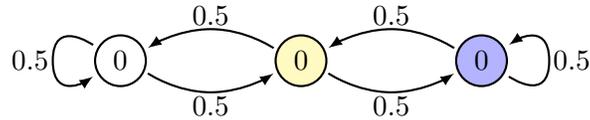
\end{example}

If one ignores the chip and focuses on the entries of the hunger vector,
the hunger game can be viewed as a greedy variant of chip-firing.
Specifically, if all entries of $\mathbf{h}$ and $P$ are rational,
we can without loss of generality assume they are all non-negative integers
(since the firing rule is invariant under affine transformation
of the hunger vector) and create a chip-configuration of the customary kind
in which the number of chips at $i$ is $\mathbf{h}_i$. 
Then the firing rule, translated to the new context,
tells us to fire the vertex $i$ that has the most chips,
with ties resolved as before.

When the state space is (countably) infinite, our hunger vectors
are infinite sequences.
We restrict ourselves to sequences that are bounded
and take on only finitely many distinct values;
this ensures that there exists at least one $i$
for which $\mathbf{h}_i$ equals $\sup_i \mathbf{h}_i$,
from which it follows that a smallest such $i$ exists.
As each vertex has only finitely many outgoing edges, 
this assumption ensures that after 
any finite number of steps in the hunger game 
there are only finitely many distinct values of hunger, 
so a vertex of highest hunger can be found
and the lowest-indexed one can be chosen, ad infinitum.

For Markov chains with absorbing states we vary the procedure slightly.
We start with a graph devoid of chips, add a chip at an initial vertex,
and then follow the rule for moving the chip described above,
with the extra stipulation that when the chip reaches 
an absorbing vertex, it gets removed from the graph.
In more detail, we define chip addition operators $E_i$ as follows:
Given an initial hunger vector $\mathbf{h}$, 
we place a chip at $v_i$ (increasing the hunger of the neighbors of $v_i$)
and add $P_i$ to $\mathbf{h}$,
and we then repeatedly move the chip to the currently hungriest vertex
(the lowest-indexed one, in the event of a tie),
simultaneously incrementing the hunger vector
by the row of $H$ corresponding to the chip's new location,
until we arrive at an absorbing vertex $v_k$,
at which point we subtract $P_k=I_k$ from the current hunger vector
and remove the chip from $v_k$.
We define $E_i(\mathbf{h})$ to be the final hunger vector,
and call $E_i$ the \textbf{chip addition operator} at $i$.
It is possible for $E_i(\mathbf{h})$ to be undefined,
in the event that the process never arrives at an absorbing state,
but we will show in \cref{lemma: finite terminate} 
that for finite absorbing chains 
the process must terminate so that $E_i$ is well-defined;
moreover, each chip addition operator preserves total hunger,
since the sum of the entries increases by 1 when $P_i$ is added,
stays the same each time a row of $H$ is added,
and decreases by 1 when $P_k$ is subtracted.

As in the previous situation, the process can be described
purely in terms of vector and matrix operations
without reference to $G$ or a chip.
Given a vector $\mathbf{h}\in\R^{n}$, add row $P_i$ to $\mathbf{h}$.
Thereafter, if $j$ is the unique value such that 
$h_{j'} < h_j$ for all $j' < j$ and $h_j \geq h_{j'}$ for all ${j'} > j$, 
add $H_j$ to $\mathbf{h}$,
unless $j$ is an absorbing state (call it $k$), in which case
subtract $P_k$ from the hunger vector and stop, 
calling the result $E_i(\mathbf{h})$.

\begin{example}\label{example: absorbing 5 walk}
Suppose we have the Markov chain represented by the graph in \cref{fig:ex1G}\,.
It has absorbing states $v_1$ and $v_5$.
\begin{figure}[htbp]
    \centering
    \begin{tikzpicture}[scale=0.6,font=\normalsize,baseline,thick]
        \foreach \x in {1,...,5} {
            \filldraw[color=black,fill=green!30,thick] (4*\x,0) circle (16pt);
            \node at (4*\x,0) {$v_{\x}$};
        }
        \draw [->,>=latex] plot [smooth,tension=5] coordinates {(8+0.5*0.7,0.866*0.7) (8+0,1.5) (8-0.5*0.7,0.866*0.7)};
        \node at (8,1.8) {0.6};
        \foreach \x in {3,...,4} {
            \draw [->,>=latex] plot [smooth,tension=1] coordinates {(4*\x-0.866*0.7,0.4*0.7) (4*\x-2,0.7) (4*\x-4+0.866*0.7,0.4*0.7)};
        }
        \foreach \x in {2,...,3} {
            \draw [->,>=latex] plot [smooth,tension=1] coordinates {(4*\x+0.866*0.7,-0.4*0.7) (4*\x+2,-0.7) (4*\x+4-0.866*0.7,-0.4*0.7)};
        }
        \draw[->,>=latex] (8-0.7,0) -- (4+0.7,0);
        \draw[->,>=latex] (16+0.7,0) -- (20-0.7,0);
        \node at (6,0.3) {0.2};
        \node at (10,1) {0.6};
        \node at (14,1) {0.2};
        \node at (10,-1) {0.2};
        \node at (14,-1) {0.4};
        \node at (18,-0.3) {0.8};
        \draw [->,>=latex] plot [smooth,tension=5] coordinates {(4-0.866*0.7,0.5*0.7) (4-1.5,0) (4-0.866*0.7,-0.5*0.7)};
        \node at (2.2,0) {1};
        \draw [->,>=latex] plot [smooth,tension=5] coordinates {(20+0.866*0.7,-0.5*0.7) (20+1.5,0) (20+0.866*0.7,0.5*0.7)};
        \node at (21.8,0) {1};
    \end{tikzpicture}
    \caption{A graph $G$ corresponding to an absorbing Markov chain.}
    \label{fig:ex1G}
\end{figure}
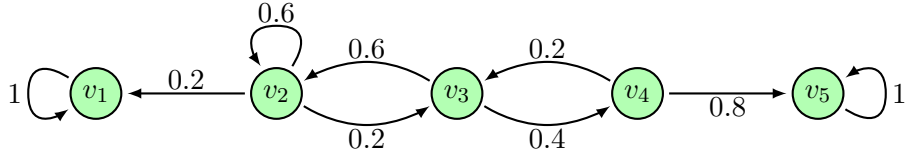
Let us compute $E_3(\mathbf{0})$.
Starting with $\mathbf{h}=\mathbf{0}$,
for our first step we add a chip to $v_3$ to yield $\mathbf{h}=[0,0.6,0,0.4,0]$,
as shown in \cref{subfig:ex1insert}\,.
After this, we follow additional steps of the hunger game process, 
firing the chip successively to $v_2$, $v_4$, and $v_5$, 
as shown in \cref{subfig:ex1fire3}\,.
Having reached an absorbing state, the final step is 
to remove the chip from $v_5$, 
resulting in $E_3(\mathbf{0})=[0.2,0.2,0.4,-0.6,-0.2]$,
as shown in \cref{subfig:ex1remove}\,.
\begin{figure}[htbp]
    \centering
    \begin{subfigure}{\textwidth}
        \centering
        \begin{tikzpicture}[scale=0.6,font=\normalsize,baseline,thick]
            \foreach \x/\xcolor in {1/white,2/white,3/white,4/white,5/white} {
                \filldraw[color=black,fill=\xcolor,thick] (4*\x,0) circle (16pt);
            }
            \foreach \x/\xtext in {1/0,2/0,3/0,4/0,5/0} {
                \node at (4*\x,0) {\xtext};
            }
            \draw [->,>=latex] plot [smooth,tension=5] coordinates {(8+0.5*0.7,0.866*0.7) (8+0,1.5) (8-0.5*0.7,0.866*0.7)};
            \node at (8,1.8) {0.6};
            \foreach \x in {3,...,4} {
                \draw [->,>=latex] plot [smooth,tension=1] coordinates {(4*\x-0.866*0.7,0.4*0.7) (4*\x-2,0.7) (4*\x-4+0.866*0.7,0.4*0.7)};
            }
            \foreach \x in {2,...,3} {
                \draw [->,>=latex] plot [smooth,tension=1] coordinates {(4*\x+0.866*0.7,-0.4*0.7) (4*\x+2,-0.7) (4*\x+4-0.866*0.7,-0.4*0.7)};
            }
            \draw[->,>=latex] (8-0.7,0) -- (4+0.7,0);
            \draw[->,>=latex] (16+0.7,0) -- (20-0.7,0);
            \node at (6,0.3) {0.2};
            \node at (10,1) {0.6};
            \node at (14,1) {0.2};
            \node at (10,-1) {0.2};
            \node at (14,-1) {0.4};
            \node at (18,-0.3) {0.8};
            \draw [->,>=latex] plot [smooth,tension=5] coordinates {(4-0.866*0.7,0.5*0.7) (4-1.5,0) (4-0.866*0.7,-0.5*0.7)};
            \node at (2.2,0) {1};
            \draw [->,>=latex] plot [smooth,tension=5] coordinates {(20+0.866*0.7,-0.5*0.7) (20+1.5,0) (20+0.866*0.7,0.5*0.7)};
            \node at (21.8,0) {1};
        \end{tikzpicture}
    \end{subfigure}
    
    \vspace{-0.5mm}
    
    \begin{subfigure}{\textwidth}
        \centering
        \begin{tikzpicture}[scale=0.6,font=\normalsize,baseline,thick]
            \foreach \x/\xcolor in {1/white,2/yellow!30,3/blue!30,4/yellow!30,5/white} {
                \filldraw[color=black,fill=\xcolor,thick] (4*\x,0) circle (16pt);
            }
            \foreach \x/\xtext in {1/0,2/0.6,3/0,4/0.4,5/0} {
                \node at (4*\x,0) {\xtext};
            }
            \draw [->,>=latex] plot [smooth,tension=5] coordinates {(8+0.5*0.7,0.866*0.7) (8+0,1.5) (8-0.5*0.7,0.866*0.7)};
            \node at (8,1.8) {0.6};
            \foreach \x in {3,...,4} {
                \draw [->,>=latex] plot [smooth,tension=1] coordinates {(4*\x-0.866*0.7,0.4*0.7) (4*\x-2,0.7) (4*\x-4+0.866*0.7,0.4*0.7)};
            }
            \foreach \x in {2,...,3} {
                \draw [->,>=latex] plot [smooth,tension=1] coordinates {(4*\x+0.866*0.7,-0.4*0.7) (4*\x+2,-0.7) (4*\x+4-0.866*0.7,-0.4*0.7)};
            }
            \draw[->,>=latex] (8-0.7,0) -- (4+0.7,0);
            \draw[->,>=latex] (16+0.7,0) -- (20-0.7,0);
            \node at (6,0.3) {0.2};
            \node at (10,1) {0.6};
            \node at (14,1) {0.2};
            \node at (10,-1) {0.2};
            \node at (14,-1) {0.4};
            \node at (18,-0.3) {0.8};
            \draw [->,>=latex] plot [smooth,tension=5] coordinates {(4-0.866*0.7,0.5*0.7) (4-1.5,0) (4-0.866*0.7,-0.5*0.7)};
            \node at (2.2,0) {1};
            \draw [->,>=latex] plot [smooth,tension=5] coordinates {(20+0.866*0.7,-0.5*0.7) (20+1.5,0) (20+0.866*0.7,0.5*0.7)};
            \node at (21.8,0) {1};
        \end{tikzpicture}
        \caption{The effect of inserting a chip at $v_3$, shown in blue. States with updated hungers are shown in yellow.}
        \label{subfig:ex1insert}
    \end{subfigure}
    \begin{subfigure}{\textwidth}
        \centering
        \begin{tikzpicture}[scale=0.6,font=\normalsize,baseline,thick]
            \foreach \x/\xcolor in {1/yellow!30,2/blue!30,3/yellow!30,4/white,5/white} {
                \filldraw[color=black,fill=\xcolor,thick] (4*\x,0) circle (16pt);
            }
            \foreach \x/\xtext in {1/0.2,2/0.2,3/0.2,4/0.4,5/0} {
                \node at (4*\x,0) {\xtext};
            }
            \draw [->,>=latex] plot [smooth,tension=5] coordinates {(8+0.5*0.7,0.866*0.7) (8+0,1.5) (8-0.5*0.7,0.866*0.7)};
            \node at (8,1.8) {0.6};
            \foreach \x in {3,...,4} {
                \draw [->,>=latex] plot [smooth,tension=1] coordinates {(4*\x-0.866*0.7,0.4*0.7) (4*\x-2,0.7) (4*\x-4+0.866*0.7,0.4*0.7)};
            }
            \foreach \x in {2,...,3} {
                \draw [->,>=latex] plot [smooth,tension=1] coordinates {(4*\x+0.866*0.7,-0.4*0.7) (4*\x+2,-0.7) (4*\x+4-0.866*0.7,-0.4*0.7)};
            }
            \draw[->,>=latex] (8-0.7,0) -- (4+0.7,0);
            \draw[->,>=latex] (16+0.7,0) -- (20-0.7,0);
            \node at (6,0.3) {0.2};
            \node at (10,1) {0.6};
            \node at (14,1) {0.2};
            \node at (10,-1) {0.2};
            \node at (14,-1) {0.4};
            \node at (18,-0.3) {0.8};
            \draw [->,>=latex] plot [smooth,tension=5] coordinates {(4-0.866*0.7,0.5*0.7) (4-1.5,0) (4-0.866*0.7,-0.5*0.7)};
            \node at (2.2,0) {1};
            \draw [->,>=latex] plot [smooth,tension=5] coordinates {(20+0.866*0.7,-0.5*0.7) (20+1.5,0) (20+0.866*0.7,0.5*0.7)};
            \node at (21.8,0) {1};
        \end{tikzpicture}
    \end{subfigure}
    
    \vspace{-0.5mm}
    
    \begin{subfigure}{\textwidth}
        \centering
        \begin{tikzpicture}[scale=0.6,font=\normalsize,baseline,thick]
            \foreach \x/\xcolor in {1/white,2/white,3/yellow!30,4/blue!30,5/yellow!30} {
                \filldraw[color=black,fill=\xcolor,thick] (4*\x,0) circle (16pt);
            }
            \foreach \x/\xtext in {1/0.2,2/0.2,3/0.4,4/-0.6,5/0.8} {
                \node at (4*\x,0) {\xtext};
            }
            \draw [->,>=latex] plot [smooth,tension=5] coordinates {(8+0.5*0.7,0.866*0.7) (8+0,1.5) (8-0.5*0.7,0.866*0.7)};
            \node at (8,1.8) {0.6};
            \foreach \x in {3,...,4} {
                \draw [->,>=latex] plot [smooth,tension=1] coordinates {(4*\x-0.866*0.7,0.4*0.7) (4*\x-2,0.7) (4*\x-4+0.866*0.7,0.4*0.7)};
            }
            \foreach \x in {2,...,3} {
                \draw [->,>=latex] plot [smooth,tension=1] coordinates {(4*\x+0.866*0.7,-0.4*0.7) (4*\x+2,-0.7) (4*\x+4-0.866*0.7,-0.4*0.7)};
            }
            \draw[->,>=latex] (8-0.7,0) -- (4+0.7,0);
            \draw[->,>=latex] (16+0.7,0) -- (20-0.7,0);
            \node at (6,0.3) {0.2};
            \node at (10,1) {0.6};
            \node at (14,1) {0.2};
            \node at (10,-1) {0.2};
            \node at (14,-1) {0.4};
            \node at (18,-0.3) {0.8};
            \draw [->,>=latex] plot [smooth,tension=5] coordinates {(4-0.866*0.7,0.5*0.7) (4-1.5,0) (4-0.866*0.7,-0.5*0.7)};
            \node at (2.2,0) {1};
            \draw [->,>=latex] plot [smooth,tension=5] coordinates {(20+0.866*0.7,-0.5*0.7) (20+1.5,0) (20+0.866*0.7,0.5*0.7)};
            \node at (21.8,0) {1};
        \end{tikzpicture}
    \end{subfigure}
    
    \vspace{-0.5mm}
    
    \begin{subfigure}{\textwidth}
        \centering
        \begin{tikzpicture}[scale=0.6,font=\normalsize,baseline,thick]
            \foreach \x/\xcolor in {1/white,2/white,3/white,4/white,5/blue!30} {
                \filldraw[color=black,fill=\xcolor,thick] (4*\x,0) circle (16pt);
            }
            \foreach \x/\xtext in {1/0.2,2/0.2,3/0.4,4/-0.6,5/0.8} {
                \node at (4*\x,0) {\xtext};
            }
            \draw [->,>=latex] plot [smooth,tension=5] coordinates {(8+0.5*0.7,0.866*0.7) (8+0,1.5) (8-0.5*0.7,0.866*0.7)};
            \node at (8,1.8) {0.6};
            \foreach \x in {3,...,4} {
                \draw [->,>=latex] plot [smooth,tension=1] coordinates {(4*\x-0.866*0.7,0.4*0.7) (4*\x-2,0.7) (4*\x-4+0.866*0.7,0.4*0.7)};
            }
            \foreach \x in {2,...,3} {
                \draw [->,>=latex] plot [smooth,tension=1] coordinates {(4*\x+0.866*0.7,-0.4*0.7) (4*\x+2,-0.7) (4*\x+4-0.866*0.7,-0.4*0.7)};
            }
            \draw[->,>=latex] (8-0.7,0) -- (4+0.7,0);
            \draw[->,>=latex] (16+0.7,0) -- (20-0.7,0);
            \node at (6,0.3) {0.2};
            \node at (10,1) {0.6};
            \node at (14,1) {0.2};
            \node at (10,-1) {0.2};
            \node at (14,-1) {0.4};
            \node at (18,-0.3) {0.8};
            \draw [->,>=latex] plot [smooth,tension=5] coordinates {(4-0.866*0.7,0.5*0.7) (4-1.5,0) (4-0.866*0.7,-0.5*0.7)};
            \node at (2.2,0) {1};
            \draw [->,>=latex] plot [smooth,tension=5] coordinates {(20+0.866*0.7,-0.5*0.7) (20+1.5,0) (20+0.866*0.7,0.5*0.7)};
            \node at (21.8,0) {1};
        \end{tikzpicture}
        \caption{$\mathbf{h}$ as the chip fires successively to $v_2$, $v_4$, and $v_5$, shown in blue. Updated hungers are shown in yellow.}
        \label{subfig:ex1fire3}
    \end{subfigure}
    \begin{subfigure}{\textwidth}
        \centering
        \begin{tikzpicture}[scale=0.6,font=\normalsize,baseline,thick]
            \foreach \x/\xcolor in {1/white,2/white,3/white,4/white,5/yellow!30} {
                \filldraw[color=black,fill=\xcolor,thick] (4*\x,0) circle (16pt);
            }
            \foreach \x/\xtext in {1/0.2,2/0.2,3/0.4,4/-0.6,5/-0.2} {
                \node at (4*\x,0) {\xtext};
            }
            \draw [->,>=latex] plot [smooth,tension=5] coordinates {(8+0.5*0.7,0.866*0.7) (8+0,1.5) (8-0.5*0.7,0.866*0.7)};
            \node at (8,1.8) {0.6};
            \foreach \x in {3,...,4} {
                \draw [->,>=latex] plot [smooth,tension=1] coordinates {(4*\x-0.866*0.7,0.4*0.7) (4*\x-2,0.7) (4*\x-4+0.866*0.7,0.4*0.7)};
            }
            \foreach \x in {2,...,3} {
                \draw [->,>=latex] plot [smooth,tension=1] coordinates {(4*\x+0.866*0.7,-0.4*0.7) (4*\x+2,-0.7) (4*\x+4-0.866*0.7,-0.4*0.7)};
            }
            \draw[->,>=latex] (8-0.7,0) -- (4+0.7,0);
            \draw[->,>=latex] (16+0.7,0) -- (20-0.7,0);
            \node at (6,0.3) {0.2};
            \node at (10,1) {0.6};
            \node at (14,1) {0.2};
            \node at (10,-1) {0.2};
            \node at (14,-1) {0.4};
            \node at (18,-0.3) {0.8};
            \draw [->,>=latex] plot [smooth,tension=5] coordinates {(4-0.866*0.7,0.5*0.7) (4-1.5,0) (4-0.866*0.7,-0.5*0.7)};
            \node at (2.2,0) {1};
            \draw [->,>=latex] plot [smooth,tension=5] coordinates {(20+0.866*0.7,-0.5*0.7) (20+1.5,0) (20+0.866*0.7,0.5*0.7)};
            \node at (21.8,0) {1};
        \end{tikzpicture}
        \caption{$\mathbf{h}$ after removing chip from $v_5$, shown in yellow.}
        \label{subfig:ex1remove}
    \end{subfigure}
    \caption{The hunger game on $G$ from \cref{fig:ex1G} after inserting a chip at $v_3$.}
    \label{fig:ex1H}
\end{figure}
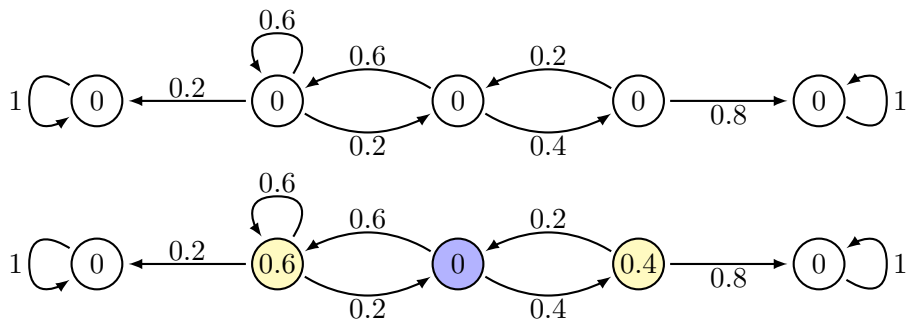
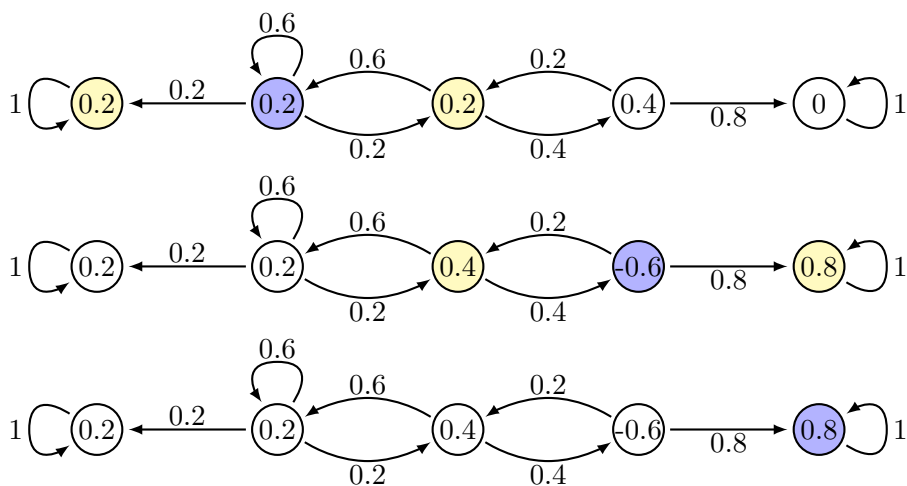
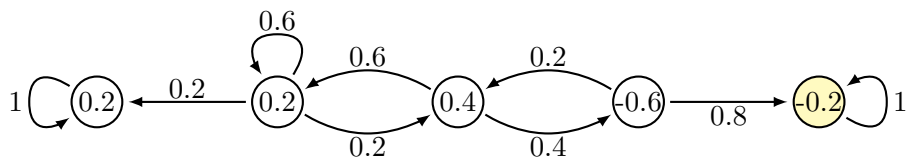

Notice that the total hunger is 0 at the start,
increases to 1 when the chip is added at $v_3$,
stays 1 as the chip moves through $G$,
and decreases to 0 when the chip is removed at the end.

\end{example}

Since increasing the hunger at every vertex by the same amount
has no effect on the dynamics of the hunger game,
when our Markov chain is finite we will often assume 
that total hunger is 0.

%% file: boundedness.tex
\section{Boundedness of hunger}\label{section: boundedness hunger}
The following lemmas will be useful in our discussion of the hunger game.
The first lemma shows that on a countable Markov chain, 
hunger stays uniformly bounded from below under the hunger game process,
including the removal and reinsertion of the chip 
if it reaches an absorbing state.

\begin{lemma}\label{lemma: hunger bounded lower}
Suppose we play the hunger game for a countable Markov chain, 
with an initial hunger vector $\mathbf{h}^{(0)}$ that is bounded below by $h_{min}\in\R$,
meaning $\h_i^{(0)} \geq h_{min}$ for all $i$.
Then hunger remains uniformly bounded below by $h_{min}-1$
under iteration of the hunger game process.
That is, suppose we have a sequence of hunger states 
$\mathbf{h}^{(1)}, \mathbf{h}^{(2)},\dots$ such that 
for each $k \geq 1$, $\mathbf{h}^{(k)}$ is reached 
from $\h^{(k-1)}$ by either firing the chip if it is not at an absorbing state
or removing and reinserting the chip if it is.
Then $\mathbf{h}^{(k)}_i \geq h_{min}-1$ holds for all $i$ and $k$.
\end{lemma}
\begin{proof}
Without loss of generality, we may take $h_{min}$ to be 0,
so that $\mathbf{h}_i^{(0)} \geq 0$ for all $i$.
Suppose our claim is false,
and let $k$ be the smallest index for which the claim fails,
so that $\h^{(k)}_i < -1$ for some $i$.
Consider the set $S$ of states whose hunger changed 
during the hunger game process from $\h^{(0)}$ to $\h^{(k)}$; 
this set must be finite, for each of the $k$ steps can only
change the hunger of finitely many states, and must contain $i$.
The total hunger of this set must be constant during this hunger game process
if we combine the process of removing a chip and inserting the next chip
into a single operation.
Let $h$ be the total hunger of $S$ throughout the process.
As the hunger of all states was originally bounded below by $0$,
we have $h \geq 0$.
Since in going from $\h^{(k-1)}$ to $\h^{(k)}$ we fired the chip to $i$,
$v_i$ must have been one of the hungriest vertices.
But in $\h^{(k-1)}$, $v_i$ must have had hunger $< 0$
(otherwise it could not have reached hunger $< -1$ in a single firing).
Since the hungers of the vertices in $S$ have nonnegative sum,
at least one of them must be nonnegative,
and hence $v_i$ could not have been one of the hungriest vertices.
This contradiction shows that hunger is bounded from below,
and more specifically that $h_{min}-1$ is a lower bound.
\end{proof}

The lemma remains true if one replaces $\geq$ by $>$
in both the hypothesis and the conclusion.

\bigskip

\begin{figure}[htbp!]
    \centering
    \begin{tikzpicture}[scale=0.5,font=\normalsize,baseline,thick]
        \filldraw[color=black,fill=green!30,thick] (0,0) circle (20pt);
        \node at (0,0) {$v_{0,0}$};
        \foreach \y in {1,...,4} {
        \foreach \x in {0,...,\y} {
            \filldraw[color=black,fill=green!30,thick] (4*\x, -4 * \y) circle (20pt);
            \node at (4*\x,-4 * \y) {$v_{\y,\x}$};
        }
        }
        \foreach \x in {0,...,4} {
            \node at (4*\x, -4*5-0.3) {$\vdots$};
        }
        \node at (4*4+2.5,-4*5+2) {$\ddots$};
        \draw [->,>=latex] plot [smooth,tension=5] coordinates {(0+0.5*0.7,0.866*0.7) (0+0,1.5) (0-0.5*0.7,0.866*0.7)};
        \node at (0,1.9) {1};
        \draw[->,>=latex] (4-0.5,-4+0.5) -- (0+0.5,0-0.5);
        \node at (2+0.3,-2+0.3) {1};
        \draw[->,>=latex] plot [smooth,tension=1] coordinates {(8,-8+0.7) (4+1,-4+1) (0+0.6,0-0.4)};
        \node at (4+1.3,-4+1.3) {1};
        \draw[->,>=latex] plot [smooth,tension=1] coordinates {(12,-12+0.7) (6+2,-6+2) (0+0.65,0-0.2)};
        \node at (6+2.3,-6+2.3) {1};
        \draw[->,>=latex] plot [smooth,tension=1] coordinates {(16,-16+0.7) (8+3,-8+3) (0+0.7,0-0)};
        \node at (8+3.3,-8+3.3) {1};
        \foreach \y in {1,...,4} {
            \draw[->,>=latex] (4*\y - 4+0.7, -4*\y) -- (4*\y -0.7, -4*\y);
            \node at (4*\y - 2,-4*\y + 0.4) {1};
        }
        \foreach \y in {2,...,4} {
            \draw[->,>=latex] (4*\y - 4*2 +0.7, -4*\y) -- (4*\y -4*1 -0.7, -4*\y);
            \node at (4*\y - 4*1 - 2, -4*\y + 0.4) {1/2};
            \draw [->,>=latex] plot [smooth,tension=1] coordinates {(4*\y-4*2+0.866*0.7,-4*\y - 0.4*0.7) (4*\y-4*1,-4*\y-0.8) (4*\y-0.866*0.7,-4*\y - 0.4*0.7)};
            \node at (4*\y-4*1,-4*\y-1.3) {1/2};
        }
        \foreach \y in {3,...,4} {
            \draw[->,>=latex] (4*\y - 4*3 +0.7, -4*\y) -- (4*\y -4*2 -0.7, -4*\y);
            \node at (4*\y - 4*2 - 2, -4*\y + 0.4) {1/3};
            \draw[->,>=latex] plot [smooth,tension=1] coordinates {(4*\y-4*3+0.866*0.7,-4*\y - 0.4*0.7) (4*\y-4*2,-4*\y-0.8) (4*\y-4*1-0.866*0.7,-4*\y - 0.4*0.7)};
            \node at (4*\y-4*2,-4*\y-1.3) {1/3};
            \draw[->,>=latex] plot [smooth,tension=1] coordinates {(4*\y-4*3+0.5,-4*\y-0.5) (4*\y-2*3,-4*\y - 2) (4*\y-0.5,-4*\y-0.5)};
            \node at (4*\y-2*3,-4*\y - 2 - 0.5) {1/3};
        }
        \foreach \y in {4,...,4} {
            \draw[->,>=latex] (4*\y - 4*4 +0.7, -4*\y) -- (4*\y -4*3 -0.7, -4*\y);
            \node at (4*\y - 4*3 - 2, -4*\y + 0.3) {1/4};
            \draw[->,>=latex] plot [smooth,tension=1] coordinates {(4*\y-4*4+0.866*0.7,-4*\y - 0.4*0.7) (4*\y-4*3,-4*\y-0.8) (4*\y-4*2-0.866*0.7,-4*\y - 0.4*0.7)};
            \node at (4*\y-4*3,-4*\y-1.3) {1/4};
            \draw[->,>=latex] plot [smooth,tension=1] coordinates {(4*\y-4*4+0.5,-4*\y-0.5) (4*\y-4*2.5,-4*\y - 2) (4*\y-4-0.5,-4*\y-0.5)};
            \node at (4*\y-2*5,-4*\y - 2 - 0.5) {1/4};
            \draw[->,>=latex] plot [smooth,tension=1] coordinates {(4*\y-4*4+0.4*0.7,-4*\y-0.866*0.7) (4*\y-4*2,-4*\y-3.1) (4*\y-4*0-0.4*0.7,-4*\y-0.866*0.7)};
            \node at (4*\y-4*2,-4*\y-3.6) {1/4};
        }
        
    \end{tikzpicture}
    \caption{A countably infinite absorbing Markov chain where hunger is not uniformly bounded under chip addition operators.}
    \label{fig:rem unbounded hunger harmonic}
\end{figure}
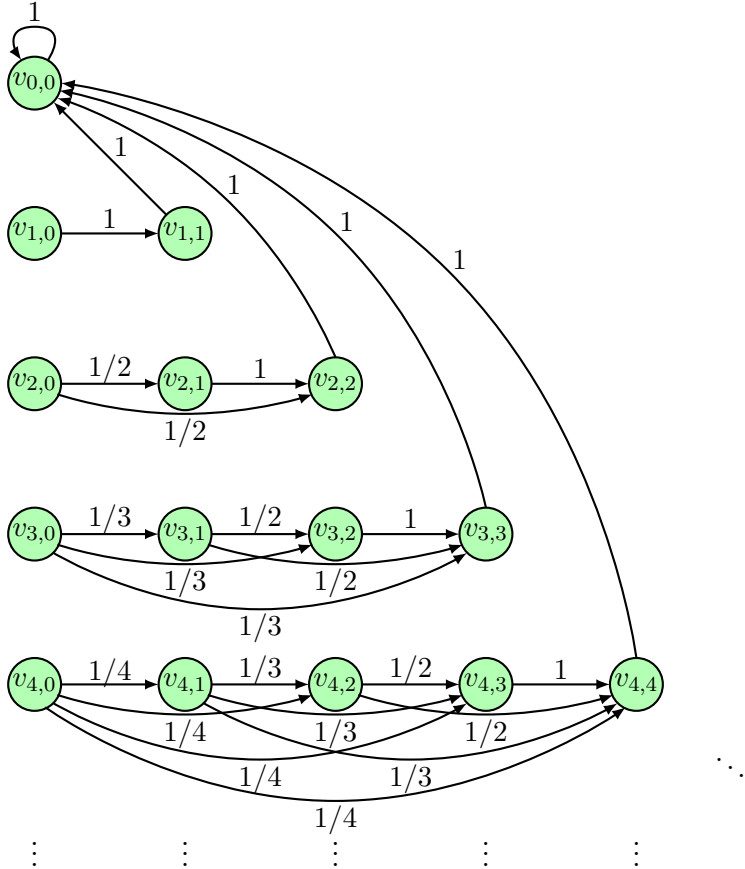

\cref{lemma: hunger bounded lower} only provides a uniform lower bound on hunger; a uniform upper bound on hunger does not exist in general for countable Markov chains.
One may construct Markov chains where depending on the location 
of chip insertion, hunger 
can grow arbitrarily large at a certain state.
As a specific example, 
consider the countably infinite Markov chain 
shown in \cref{fig:rem unbounded hunger harmonic}\,.
We have states $v_{i,j}$ for all integers $i \geq j \geq 0$, 
where $v_{0,0}$ is an absorbing state.
A walker at vertex $v_{i,j}$ for $i > j$ moves with equal probability 
to any vertex $v_{i,k}$ with $k>j$, 
in other words with probability $\frac{1}{i-j}$ to $v_{i,k}$ 
for all $j+1 \leq k \leq i$, 
and a walker at $v_{i,i}$ moves with probability 1 
to the absorbing vertex $v_{0,0}$.
Starting with $\mathbf{h}=\mathbf{0}$, 
inserting a chip at $v_{n,0}$ for $n>0$ 
causes the chip to move from $v_{n,0}$ to $v_{n,n}$, 
increasing the second coordinate by 1 each step 
(i.e. moving one position to the right at each step).
When the chip fires to $v_{n,n-1}$, 
state $v_{n,n}$ has hunger $1+\frac{1}{2}+\frac{1}{3}+\cdots+\frac{1}{n}$, 
the $n$th partial sum of the harmonic series. 
As the harmonic series diverges, by picking sufficiently large $n$, 
hunger can become arbitrarily large at $v_{n,n}$, 
and is thus not uniformly bounded under the hunger game process,
when chip removal and insertion are allowed.

This example used a countably infinite state space, which is necessary for this unbounded behavior to occur.
As shown in the following lemma, for a finite Markov chain, 
hunger stays uniformly bounded from both sides under the hunger game process,
including the removal and reinsertion of the chip 
if it reaches an absorbing state.

\begin{lemma}\label{lemma: hunger bounded}
For a given hunger state $\h^{(0)}$ on a finite Markov chain, 
hunger remains uniformly bounded under the hunger game process.
In other words, there exist $a,b \in \R$ such that for any sequence
of hunger states $\h^{(1)},\h^{(2)},\dots$, where for each $k\geq 1$,
$\h^{(k)}$ is reached from $\h^{(k-1)}$ by either
firing the chip if it is not at an absorbing state
or removing and reinserting the chip if it is,
the inequality $a \leq \h^{(k)}_i \leq b$ holds for all $i$ and $k$.
\end{lemma}
\begin{proof}
The existence of a lower bound $a$ follows directly 
from \cref{lemma: hunger bounded lower}\,, 
as any hunger state on a finite Markov chain must be bounded.

As hunger is bounded from below by $a$ 
and total hunger $h$ is constant during the hunger game process, 
when we conceive of removing a chip and inserting the next chip 
as occurring simultaneously, 
the hunger of any single vertex cannot exceed $h-(n-1)a$, 
and thus hunger is bounded from above and $b$ exists.
\end{proof}

\begin{remark}\label{remark: finite orbit}
When all transition probabilities are rational
and the Markov chain has no absorbing states,
so that we never remove and reinsert the chip---which
introduces choice into the hunger game---we claim \cref{lemma: hunger bounded}
implies that the hunger game is eventually periodic.
Let $d$ be the least common multiple
of the denominators of all the transition probabilities.
For any $\h^{(0)}$, as each entry is uniformly bounded and may only
change by a multiple of $\frac{1}{d}$, there are only finitely
many values that $\h^{(k)}$ may take on as $k$ varies.
This implies that eventually we must have $\h^{(j)} = \h^{(i)}$ with $j > i$,
so that the hunger game has become periodic with period dividing $j-i$.
If we define $\v$ as the vector that counts
the number of visits to each state of the chain from time $i+1$ to time $j$,
we have $\v H = {\bf 0}$, implying that $\v$ is a stationary vector
and that $\frac{1}{j-i}\v$ is the stationary probability measure for the Markov chain.
Theorem \ref{theorem: irreducible fidelity converge stationary}
extends this claim, in a suitable asymptotic sense,
to situations in which the transition probabilities
are not all rational.
\end{remark}

%% file: termination.tex
\section{Termination of stabilization processes}\label{section: termination}
An \textbf{absorbing Markov chain} is a (not necessarily finite) Markov chain 
where each vertex has a path of finite length to an absorbing state.
For example, the infinite Markov chain 
shown in \cref{fig:rem unbounded hunger harmonic}
is absorbing, as a walker at any vertex $v_{i,j}$ can move 
to absorbing vertex $v_{0,0}$ within 2 moves.
We claim that in a finite absorbing Markov chain, 
the chip addition process at any vertex $v$ in the corresponding graph $G$ 
is guaranteed to terminate, regardless of the initial hunger state, 
and thus the chip addition operators $E_i$ are well-defined
for finite absorbing Markov chains.
This follows immediately
from the finiteness of $G$ and the following lemma:

\begin{lemma}\label{lemma: finite terminate}
For a finite absorbing Markov chain, each vertex $v$ gets visited
only finitely often in the hunger game
before the chip reaches an absorbing vertex.
\end{lemma}

\begin{proof}
The assumption that the Markov chain is absorbing
guarantees that for each $v$
there is a finite path $w_0,w_1,\dots,w_m$ 
such that $w_0$ is $v$, $w_m$ is an absorbing vertex,
and the transition probability from $w_i$ to $w_{i+1}$
is positive for all $0 \leq i \leq m-1$.
We prove the claim by induction on $m$.
The case $m=0$ is trivial,
as once the chip visits an absorbing vertex it is removed.
For a less trivial case, consider $m=1$.
For sake of contradiction, suppose $v=w_0$ is visited infinitely often.
Each time $w_0$ is visited, the hunger of $w_1$ increases by a fixed amount.
As $w_1$ is an absorbing vertex, we cannot visit it, as otherwise we would only visit $w_0$ finitely often before reaching an absorbing vertex.
This implies the hunger of $w_1$ increases unboundedly due to $w_0$ being visited infinitely often, contradicting the existence of an upper bound on hunger demonstrated in \cref{lemma: hunger bounded}.
Hence $w_0$ is visited finitely often, proving the case $m=1$.

The inductive step follows similar reasoning to the $m=1$ argument.
Suppose now that the claim is true for $m-1$,
so that $w_1$ is visited only finitely often.
Assume for sake of contradiction that $w_0$ is visited infinitely often, 
where each visit increases the hunger of $w_1$ by a fixed amount.
After $w_1$ is visited for the last time,
its hunger must grow without bound due to $w_0$
being visited infinitely often;
this contradicts \cref{lemma: hunger bounded}\,.
Hence the claim holds for all vertices.
\end{proof}

Since $G$ has only finitely many vertices,
each of which can fire only finitely many times before the chip is absorbed,
the chip must eventually be absorbed, as claimed.

A weakened version of this lemma applies to countable absorbing Markov chains.

\begin{proposition}\label{proposition: countable terminate}
On a countable absorbing Markov chain, if we add a chip
and then perform repeated firing, then either the chip eventually gets absorbed
or the chip is not confined within any finite subset of the state space.
\end{proposition}
\begin{proof}
If the Markov chain is finite, the result immediately follows 
from \cref{lemma: finite terminate}\,,
so assume the state space is countably infinite.
Suppose there is a finite subset $U$ of the state space $V$
such that the chip is always in $U$.
Consider the set $S$ given by
\begin{align*}
    S = U \cup \{ v \in V \mid \exists\, i\in U : P_{ij}>0\}.
\end{align*}
Because $U$ is finite and each vertex has finitely many outgoing edges, 
$S$ is also finite.
Crucially, as the chip is bounded within $U$, 
the only states whose hunger can change during the chip addition process 
are the states in $S$.
Hence, we may equivalently view the hunger game as acting upon 
the induced subgraph of $G$ formed from using $S$ as the vertex set, 
and from \cref{lemma: finite terminate} the result follows.
\end{proof}

\begin{remark}\label{remark: non-terminating infinite chip addition}
For infinite absorbing Markov chains,
it is possible for the chip to wander off
without being confined to any finite subset of the state space.
As a simple example, based on the goldbug system 
introduced by the second author of this article
as described in Kleber \cite{kleber2005goldbug}, 
take $\N\cup\{-1\}$ as the state space 
where states $-1$ and 0 are absorbing, and for $i \geq 1$, 
a walker at vertex $i$ 
moves to vertex $i-2$ and $i+1$ each with probability $\frac{1}{2}$, 
as shown in \cref{fig: goldbug system}\,.
Using the hunger state $\h$ where $\h_{-1} = \h_0 = \h_1 = -\frac{1}{2}$ 
and $\h_i = 0$ for $i > 1$, as shown in \cref{subfig: goldbug init}\,,
a chip inserted at state 1, as shown in \cref{subfig: goldbug insert}\,,
will move rightwards one state at a time to infinity;
the first few steps of the process 
are shown in \cref{subfig: goldbug fire 5}\,.
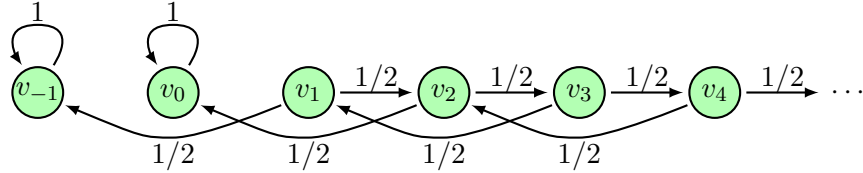
\begin{figure}[htbp]
    \centering
    \begin{tikzpicture}[scale=0.6,font=\normalsize,baseline,thick]
        \foreach \x in {-1,...,4} {
            \filldraw[color=black,fill=green!30,thick] (3*\x + 3, 0) circle (16pt);
            \node at (3*\x + 3, 0) {$v_{\x}$};
        }
        \node at (3*6,0) {$\cdots$};
        \foreach \x in {0,...,1} {
            \draw [->,>=latex] plot [smooth,tension=5] coordinates {(3*\x+0.5*0.7,0.866*0.7) (3*\x+0,1.5) (3*\x-0.5*0.7,0.866*0.7)};
            \node at (3*\x,1.8) {1};
        }
        \foreach \x in {2,...,5} {
            \draw[->,>=latex] (3*\x +0.7, 0) -- (3*\x + 3 -0.7, 0);
            \node at (3*\x +1.5,0.3) {1/2};
        }
        \foreach \x in {2,...,5} {
            \draw [<-,>=latex] plot [smooth,tension=0.5] coordinates {(3*\x-2*3+0.866*0.7,-0.4*0.7) (3*\x-1*3,-1) (3*\x-0.866*0.7,-0.4*0.7)};
            \node at (3*\x-1*3,-1.4) {1/2};
        }
    \end{tikzpicture}
    \caption{The goldbug system.}
    \label{fig: goldbug system}
\end{figure}

\begin{figure}[htbp]
    \centering
    \begin{subfigure}{\textwidth}
    \centering
    \begin{tikzpicture}[scale=0.6,font=\normalsize,baseline,thick]
        \foreach \x in {0,...,5} {
            \filldraw[color=black,fill=white,thick] (3*\x, 0) circle (16pt);
        }
        \foreach \x in {0,...,2} {
            \node at (3*\x,0) {$-\frac{1}{2}$};
        }
        \foreach \x in {3,...,5} {
            \node at (3*\x,0) {0};
        }
        \node at (3*6,0) {$\cdots$};
        \foreach \x in {0,...,1} {
            \draw [->,>=latex] plot [smooth,tension=5] coordinates {(3*\x+0.5*0.7,0.866*0.7) (3*\x+0,1.5) (3*\x-0.5*0.7,0.866*0.7)};
            \node at (3*\x,1.8) {1};
        }
        \foreach \x in {2,...,5} {
            \draw[->,>=latex] (3*\x +0.7, 0) -- (3*\x + 3 -0.7, 0);
            \node at (3*\x +1.5,0.3) {1/2};
        }
        \foreach \x in {2,...,5} {
            \draw [<-,>=latex] plot [smooth,tension=0.5] coordinates {(3*\x-2*3+0.866*0.7,-0.4*0.7) (3*\x-1*3,-1) (3*\x-0.866*0.7,-0.4*0.7)};
            \node at (3*\x-1*3,-1.4) {1/2};
        }
    \end{tikzpicture}
    \caption{The initial hunger state $\h$.}
    \label{subfig: goldbug init}
    \end{subfigure}
    \begin{subfigure}{\textwidth}
    \centering
    \begin{tikzpicture}[scale=0.6,font=\normalsize,baseline,thick]
        \foreach \x in {0,...,5} {
            \filldraw[color=black,fill=white,thick] (3*\x, 0) circle (16pt);
        }
        \filldraw[color=black,fill=blue!30,thick] (3*2,0) circle (16pt);
        \filldraw[color=black,fill=yellow!30,thick] (3*0,0) circle (16pt);
        \filldraw[color=black,fill=yellow!30,thick] (3*3,0) circle (16pt);
        \node at (3*0,0) {$0$};
        \node at (3*1,0) {$-\frac{1}{2}$};
        \node at (3*2,0) {$-\frac{1}{2}$};
        \node at (3*3,0) {$\frac{1}{2}$};
        \node at (3*4,0) {$0$};
        \node at (3*5,0) {$0$};
        \node at (3*6,0) {$\cdots$};
        \foreach \x in {0,...,1} {
            \draw [->,>=latex] plot [smooth,tension=5] coordinates {(3*\x+0.5*0.7,0.866*0.7) (3*\x+0,1.5) (3*\x-0.5*0.7,0.866*0.7)};
            \node at (3*\x,1.8) {1};
        }
        \foreach \x in {2,...,5} {
            \draw[->,>=latex] (3*\x +0.7, 0) -- (3*\x + 3 -0.7, 0);
            \node at (3*\x +1.5,0.3) {1/2};
        }
        \foreach \x in {2,...,5} {
            \draw [<-,>=latex] plot [smooth,tension=0.5] coordinates {(3*\x-2*3+0.866*0.7,-0.4*0.7) (3*\x-1*3,-1) (3*\x-0.866*0.7,-0.4*0.7)};
            \node at (3*\x-1*3,-1.4) {1/2};
        }
    \end{tikzpicture}
    \caption{$\h$ after inserting a chip at $v_1$, shown in blue. States with updated hungers are shown in yellow.}
    \label{subfig: goldbug insert}
    \end{subfigure}
    \begin{subfigure}{\textwidth}
    \centering
    \begin{tikzpicture}[scale=0.6,font=\normalsize,baseline,thick]
        \foreach \x in {0,...,5} {
            \filldraw[color=black,fill=white,thick] (3*\x, 0) circle (16pt);
        }
        \filldraw[color=black,fill=blue!30,thick] (3*3,0) circle (16pt);
        \filldraw[color=black,fill=yellow!30,thick] (3*1,0) circle (16pt);
        \filldraw[color=black,fill=yellow!30,thick] (3*4,0) circle (16pt);
        \node at (3*0,0) {$0$};
        \node at (3*1,0) {$0$};
        \node at (3*2,0) {$-\frac{1}{2}$};
        \node at (3*3,0) {$-\frac{1}{2}$};
        \node at (3*4,0) {$\frac{1}{2}$};
        \node at (3*5,0) {$0$};
        \node at (3*6,0) {$\cdots$};
        \foreach \x in {0,...,1} {
            \draw [->,>=latex] plot [smooth,tension=5] coordinates {(3*\x+0.5*0.7,0.866*0.7) (3*\x+0,1.5) (3*\x-0.5*0.7,0.866*0.7)};
            \node at (3*\x,1.8) {1};
        }
        \foreach \x in {2,...,5} {
            \draw[->,>=latex] (3*\x +0.7, 0) -- (3*\x + 3 -0.7, 0);
            \node at (3*\x +1.5,0.3) {1/2};
        }
        \foreach \x in {2,...,5} {
            \draw [<-,>=latex] plot [smooth,tension=0.5] coordinates {(3*\x-2*3+0.866*0.7,-0.4*0.7) (3*\x-1*3,-1) (3*\x-0.866*0.7,-0.4*0.7)};
            \node at (3*\x-1*3,-1.4) {1/2};
        }
    \end{tikzpicture}
    \end{subfigure}
    \begin{subfigure}{\textwidth}
    \centering
    \begin{tikzpicture}[scale=0.6,font=\normalsize,baseline,thick]
        \foreach \x in {0,...,5} {
            \filldraw[color=black,fill=white,thick] (3*\x, 0) circle (16pt);
        }
        \filldraw[color=black,fill=blue!30,thick] (3*4,0) circle (16pt);
        \filldraw[color=black,fill=yellow!30,thick] (3*2,0) circle (16pt);
        \filldraw[color=black,fill=yellow!30,thick] (3*5,0) circle (16pt);
        \node at (3*0,0) {$0$};
        \node at (3*1,0) {$0$};
        \node at (3*2,0) {$0$};
        \node at (3*3,0) {$-\frac{1}{2}$};
        \node at (3*4,0) {$-\frac{1}{2}$};
        \node at (3*5,0) {$\frac{1}{2}$};
        \node at (3*6,0) {$\cdots$};
        \foreach \x in {0,...,1} {
            \draw [->,>=latex] plot [smooth,tension=5] coordinates {(3*\x+0.5*0.7,0.866*0.7) (3*\x+0,1.5) (3*\x-0.5*0.7,0.866*0.7)};
            \node at (3*\x,1.8) {1};
        }
        \foreach \x in {2,...,5} {
            \draw[->,>=latex] (3*\x +0.7, 0) -- (3*\x + 3 -0.7, 0);
            \node at (3*\x +1.5,0.3) {1/2};
        }
        \foreach \x in {2,...,5} {
            \draw [<-,>=latex] plot [smooth,tension=0.5] coordinates {(3*\x-2*3+0.866*0.7,-0.4*0.7) (3*\x-1*3,-1) (3*\x-0.866*0.7,-0.4*0.7)};
            \node at (3*\x-1*3,-1.4) {1/2};
        }
    \end{tikzpicture}
    \end{subfigure}
    \begin{subfigure}{\textwidth}
    \centering
    \begin{tikzpicture}[scale=0.6,font=\normalsize,baseline,thick]
        \foreach \x in {0,...,5} {
            \filldraw[color=black,fill=white,thick] (3*\x, 0) circle (16pt);
        }
        \filldraw[color=black,fill=blue!30,thick] (3*5,0) circle (16pt);
        \filldraw[color=black,fill=yellow!30,thick] (3*3,0) circle (16pt);
        \node at (3*0,0) {$0$};
        \node at (3*1,0) {$0$};
        \node at (3*2,0) {$0$};
        \node at (3*3,0) {$0$};
        \node at (3*4,0) {$-\frac{1}{2}$};
        \node at (3*5,0) {$-\frac{1}{2}$};
        \node at (3*6,0) {$\cdots$};
        \foreach \x in {0,...,1} {
            \draw [->,>=latex] plot [smooth,tension=5] coordinates {(3*\x+0.5*0.7,0.866*0.7) (3*\x+0,1.5) (3*\x-0.5*0.7,0.866*0.7)};
            \node at (3*\x,1.8) {1};
        }
        \foreach \x in {2,...,5} {
            \draw[->,>=latex] (3*\x +0.7, 0) -- (3*\x + 3 -0.7, 0);
            \node at (3*\x +1.5,0.3) {1/2};
        }
        \foreach \x in {2,...,5} {
            \draw [<-,>=latex] plot [smooth,tension=0.5] coordinates {(3*\x-2*3+0.866*0.7,-0.4*0.7) (3*\x-1*3,-1) (3*\x-0.866*0.7,-0.4*0.7)};
            \node at (3*\x-1*3,-1.4) {1/2};
        }
    \end{tikzpicture}
    \caption{$\h$ as the chip fires successively to $v_2$, $v_3$, and $v_4$, shown in blue. Updated hungers are shown in yellow.}
    \label{subfig: goldbug fire 5}
    \end{subfigure}
    \caption{A hunger state on the goldbug system where a chip inserted at $v_1$ goes to infinity.}
    \label{fig:rem non-terminating goldbug}
\end{figure}
\end{remark}

%% file: stationary.tex
\section{Stationary distributions}\label{section: stationary distributions}
We now consider the hunger game process for nonabsorbing Markov chains, 
so that chips are never inserted or removed.

Say that a vector $\v$ is stationary under $P$ if $\v P = \v$,
or equivalently $\v H = \v(P-I) = {\bf 0}$.
Let $E$ be the space of vectors that are stationary under $P$,
or equivalently the nullspace of $H$.
When the Markov chain is irreducible,
so that there is a unique stationary distribution $\ppi$,
$E$ is the 1-dimensional subspace of $\R^n$
consisting of multiples of $[\ppi(1),\dots,\ppi(n)]$.

Define the \textbf{firing vector} $\v$ 
after $N$ steps of a hunger game process 
to be the vector whose $i$th component $\v_i$ is 
the number of times $H_i$ was added to $\h$, 
i.e.\ the number of times the chip fired to vertex $v_i$.
The following theorem demonstrates that 
the normalized firing vector $\frac{\v}{N}$ 
approximates the unique stationary distribution $\ppi$ 
of an irreducible finite Markov chain 
within a distance proportional to $N^{-1}$.

\begin{theorem}\label{theorem: irreducible fidelity converge stationary}
Given a hunger game process on an irreducible finite Markov chain, 
let $\v^{(N)}$ be the firing vector after $N$ steps.
Then the sequence of normalized firing vectors 
$\left\{\frac{1}{N}\v^{(N)}\right\}$ converges to 
the unique stationary distribution $\ppi$, 
where there exists a constant $C$ such that for all $N$, 
the normalized firing vector is within distance $\frac{C}{N}$ of $\ppi$ 
in the $L^1$ metric.
\end{theorem}
\begin{proof}
Since the Markov chain is irreducible, 0 is a simple eigenvalue of $H$.
Let $c$ be the maximum of $1/|\lambda|$
where $\lambda$ ranges over the nonzero eigenvalues of $H$.
For $d \in \R$, let $U_d = \{\x \mid \x_1 + \cdots + \x_n = d\}$, 
where there are $n$ states in the Markov chain.
The fact that 0 is a simple eigenvalue of $H$
and that every other eigenvalue has norm at least $1/c$
implies that the restriction of $H$
to an affine map from $U_1$ to $U_0$ is invertible,
and that if two points $\p^{(1)}, \p^{(2)}\in U_0$ 
are within distance $\varepsilon$, then 
their two preimages $\x^{(1)},\x^{(2)}\in U_1$
are within distance $c\varepsilon$ of each other.

It follows that if we have a point $\p\in U_0$ 
within distance $\varepsilon$ of $\mathbf{0}$, 
the preimages $\v,\ppi\in U_1$ of $\p$ and $\mathbf{0}$, respectively, 
are within distance $c\varepsilon$ of each other.
As a result, $\v$ is within distance $c\varepsilon$ 
of the stationary distribution $\ppi$.

By \cref{lemma: hunger bounded}\,,
during any hunger game process hunger remains bounded, 
so the change in hunger after $N$ steps, given by $\v^{(N)} H$, 
is within a bounded distance, say $b$, of $\mathbf{0}$.
This implies $\left(\frac{1}{N}\v^{(N)}\right)H$ is 
within distance $\frac{b}{N}$ of $\mathbf{0}$, 
which implies the normalized firing vector $\frac{1}{N}\v^{(N)}$ 
is within distance $\frac{bc}{N}$ of a stationary distribution $\ppi$.
Setting $C=bc$ yields our desired result.
\end{proof}
\begin{remark}\label{remark: fidelity converge stationary single absorbing}
As the only time irreducibility was assumed 
in \cref{theorem: irreducible fidelity converge stationary} was 
when we claimed the Markov chain had a unique stationary distribution $\ppi$,
the result holds for any finite Markov chain with a unique stationary distribution.
Namely, it also works for absorbing Markov chains 
with a single absorbing state $v_k$, where $\ppi$ is the vector 
with a 1 corresponding to state $v_k$ and 0s elsewhere.
\end{remark}

When a unique stationary distribution exists, 
the stationary distribution represents, in the long run, 
the probability distribution of the chain being at each particular state, 
or the occupation frequency distribution.
As a result, \cref{theorem: irreducible fidelity converge stationary} 
illustrates that the normalized firing vector, 
which counts how many times the chip arrives at each state in the hunger game, 
approximates the occupation frequency distribution 
within a discrepancy proportional to $N^{-1}$, 
better than the $N^{-1/2}$ discrepancy expected 
from repeated samplings of the corresponding random process.

\begin{example}\label{example: stationary distribution}
The Markov chain described in \cref{example: 3 hunger} is irreducible, 
and has unique stationary distribution 
$\ppi=\left[\frac{1}{3},\frac{1}{3},\frac{1}{3}\right]$, 
as shown in \cref{fig:ex stationary}\,.
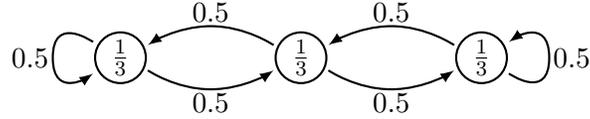
\begin{figure}[htbp]
    \centering
    \begin{tikzpicture}[scale=0.6,font=\normalsize,baseline,thick]
        \foreach \x in {1,...,3} {
            \filldraw[color=black,fill=white,thick] (4*\x,0) circle (16pt);
            \node at (4*\x,0) {$\frac{1}{3}$};
        }
        \foreach \x in {2,3} {
            \draw [->,>=latex] plot [smooth,tension=1] coordinates {(4*\x-0.866*0.7,0.4*0.7) (4*\x-2,0.7) (4*\x-4+0.866*0.7,0.4*0.7)};
        }
        \foreach \x in {1,2} {
            \draw [->,>=latex] plot [smooth,tension=1] coordinates {(4*\x+0.866*0.7,-0.4*0.7) (4*\x+2,-0.7) (4*\x+4-0.866*0.7,-0.4*0.7)};
        }
        \node at (6,1) {0.5};
        \node at (10,1) {0.5};
        \node at (6,-1) {0.5};
        \node at (10,-1) {0.5};
        \draw [->,>=latex] plot [smooth,tension=5] coordinates {(4-0.866*0.7,0.5*0.7) (4-1.5,0) (4-0.866*0.7,-0.5*0.7)};
        \node at (2,0) {0.5};
        \draw [->,>=latex] plot [smooth,tension=5] coordinates {(12+0.866*0.7,-0.5*0.7) (12+1.5,0) (12+0.866*0.7,0.5*0.7)};
        \node at (14,0) {0.5};
    \end{tikzpicture}
    \caption{The unique stationary distribution $\ppi$ of a doubly-reflecting random walk}
    \label{fig:ex stationary}
\end{figure}
From \cref{example: 3 hunger} we know that 
the states are visited periodically in the order 1, 2, and 3, 
so the firing vector starting at $\mathbf{h}=\mathbf{0}$ is given by
\begin{align*}
    \mathbf{v}^{(N)} = \left[\left\lfloor\dfrac{N+2}{3}\right\rfloor,\left\lfloor\dfrac{N+1}{3}\right\rfloor,\bigg\lfloor\dfrac{N}{3}\bigg\rfloor\right],
\end{align*}
or equivalently
\begin{align*}
    \mathbf{v}^{(3M)} &= \left[M,M,M\right] \\
    \mathbf{v}^{(3M+1)} &= \left[M+1,M,M\right] \\
    \mathbf{v}^{(3M+2)} &= \left[M+1,M+1,M\right].
\end{align*}
\cref{theorem: irreducible fidelity converge stationary} guarantees 
the existence of some constant $C$ such that 
the normalized firing vector $\frac{1}{N}\mathbf{v}^{(N)}$ 
is within $\frac{C}{N}$ of $\ppi$ for all $N$.
In fact, we will show $C=\frac{4}{3}$ is the minimal such constant.
Then for $N=3M$, the normalized firing vector is 
$\left[\frac{1}{3},\frac{1}{3},\frac{1}{3}\right]$, so there is no discrepancy,
while for $N=3M+1$ the normalized firing vector is 
$\left[\frac{M+1}{3M+1},\frac{M}{3M+1},\frac{M}{3M+1}\right]$, 
whose discrepancy from $\ppi=\left[\frac{1}{3},\frac{1}{3},\frac{1}{3}\right]$ 
is $\frac{4/3}{3M+1}=\frac{C}{N}$.
Similarly for $N=3M+2$.

In general, however, unlike in this example, 
the sequence of states visited need not be periodic, 
as the transition probabilities can be irrational.
\end{example}

%% file: hitting.tex
\section{Hitting Probabilities and Absorption Times}\label{section: hitting probabilities and hitting times}
Let $X_0, X_1, \dots$ be a Markov chain on a finite state space $V$, 
and let $\P_v$ denote the law of the Markov chain initialized at state $v$, 
or in other words $X_0 = v$.
Define $T_u$ to be the \textbf{hitting time} of vertex $u$ by the Markov chain,
given by
\begin{align*}
    T_u = \min\{t \geq 0: X_t = u\},
\end{align*}
where if no such $t$ exists then $T_u = \infty$.

For an absorbing Markov chain, 
let the set of absorbing states be $U \subseteq V$.
If our starting vertex $v$ belongs to $U$, then the behavior of the chain is trivial,
so we will assume henceforth that $v \not\in U$.
We define the \textbf{hitting probability} of $u \in U$ by
\begin{align*}
    h_u(v) = \P_v (\mbox{$X_t = u$ for all $t$ sufficiently large}),
\end{align*}
where $\P_v$ denotes the law of the Markov chain initiated at $v$.
The hitting probability is equivalently
\begin{align*}
    h_u(v) = \P_v (\forall u' \in U,\,\, T_u \leq T_{u'} ),
\end{align*}
which will give us the flexibility later on to modify the Markov chain so that the states in $U$ are no longer absorbing.
As an absorbing Markov chain reaches an absorbing state with probability 1,
for any $v$ we have
\begin{align*}
    \sum_{u \in U} h_u(v) = 1.
\end{align*}

In order to measure hitting probabilities $h_{u}(v)$ 
through the hunger game process, 
we consider a modification of the Markov chain
in which all absorbing states go to $v$ with probability 1.
In other words, for each $u\in U$ we replace $P_u$, the row of $P$ corresponding to state $u$, with $e_v$, the unit vector corresponding to $v$.
Crucially, this modification does not change 
the hitting probabilities $h_{u}(v)$
as given by our second definition.
Additionally, we remove any states in $V$ 
that cannot be reached from $v$ with positive probability;
this does not alter the hitting probabilities.
We will refer to this modified Markov chain as the 
\textbf{rerouted Markov chain} associated with $v$.
Due to the removal of unreachable states, a rerouted Markov chain of an absorbing Markov chain is always irreducible.

\begin{example}\label{example: rerouted MC}
Consider the absorbing Markov chain given in \cref{subfig: ex rerouted base}\,, 
which has absorbing states $v_1$ and $v_4$.
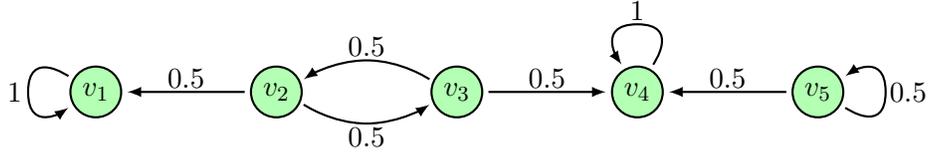
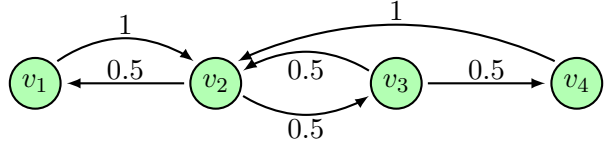
\begin{figure}[htbp]
    \centering
    \begin{subfigure}{\textwidth}
        \centering
        \begin{tikzpicture}[scale=0.6,font=\normalsize,baseline,thick]
            \foreach \x in {1,...,5} {
                \filldraw[color=black,fill=green!30,thick] (4*\x,0) circle (16pt);
                \node at (4*\x,0) {$v_{\x}$};
            }
            \draw [->,>=latex] plot [smooth,tension=5] coordinates {(4-0.866*0.7,0.5*0.7) (4-1.5,0) (4-0.866*0.7,-0.5*0.7)};
            \node at (2.2,0) {1};
            \draw[->,>=latex] (4*2-0.7,0) -- (4*1+0.7,0);
            \node at (4*1+2,0.3) {0.5};
            \draw [->,>=latex] plot [smooth,tension=1] coordinates {(4*3-0.866*0.7,0.4*0.7) (4*3-2,0.7) (4*3-4+0.866*0.7,0.4*0.7)};
            \node at (4*3-2,1) {0.5};
            \draw [->,>=latex] plot [smooth,tension=1] coordinates {(4*2+0.866*0.7,-0.4*0.7) (4*2+2,-0.7) (4*2+4-0.866*0.7,-0.4*0.7)};
            \node at (4*2+2,-1) {0.5};
            \draw[->,>=latex] (4*3+0.7,0) -- (4*4-0.7,0);
            \node at (4*3+2,0.3) {0.5};
            \draw [->,>=latex] plot [smooth,tension=5] coordinates {(4*4+0.5*0.7,0.866*0.7) (4*4+0,1.5) (4*4-0.5*0.7,0.866*0.7)};
            \node at (4*4,1.8) {1};
            \draw[->,>=latex] (4*5-0.7,0) -- (4*4+0.7,0);
            \node at (4*5-2,0.3) {0.5};
            \draw [->,>=latex] plot [smooth,tension=5] coordinates {(4*5+0.866*0.7,-0.5*0.7) (4*5+1.5,0) (4*5+0.866*0.7,0.5*0.7)};
            \node at (22,0) {0.5};
        \end{tikzpicture}
        \caption{An absorbing Markov chain.}
        \label{subfig: ex rerouted base}
    \end{subfigure}
    \begin{subfigure}{\textwidth}
        \centering
        \begin{tikzpicture}[scale=0.6,font=\normalsize,baseline,thick]
            \foreach \x in {1,...,4} {
                \filldraw[color=black,fill=green!30,thick] (4*\x,0) circle (16pt);
                \node at (4*\x,0) {$v_{\x}$};
            }
            \draw [->,>=latex] plot [smooth,tension=1] coordinates {(4*1+0.7*0.7,0.7*0.7) (4*2-2,1) (4*2-0.7*0.7,0.7*0.7)};
            \node at (4*1+2,1.3) {1};
            \draw[->,>=latex] (4*2-0.7,0) -- (4*1+0.7,0);
            \node at (4*1+2,0.3) {0.5};
            \draw [->,>=latex] plot [smooth,tension=1] coordinates {(4*3-0.866*0.7,0.4*0.7) (4*3-2,0.7) (4*3-4+0.866*0.7,0.4*0.7)};
            \node at (4*3-2,0.3) {0.5};
            \draw [->,>=latex] plot [smooth,tension=1] coordinates {(4*2+0.866*0.7,-0.4*0.7) (4*2+2,-0.7) (4*2+4-0.866*0.7,-0.4*0.7)};
            \node at (4*2+2,-1) {0.5};
            \draw[->,>=latex] (4*3+0.7,0) -- (4*4-0.7,0);
            \node at (4*3+2,0.3) {0.5};
            \draw [->,>=latex] plot [smooth,tension=1] coordinates {(4*4-0.7*0.7,0.7*0.7) (4*3,1.3) (4*2+0.7*0.7,0.7*0.7)};
            \node at (4*3,1.6) {1};
        \end{tikzpicture}
        \caption{The rerouted Markov chain associated with $v_2$.}
        \label{subfig: ex rerouted 2}
    \end{subfigure}
    \caption{A rerouted Markov chain.}
    \label{fig: ex rerouted}
\end{figure}
The rerouted Markov chain associated with non-absorbing state $v_2$ 
is shown in \cref{subfig: ex rerouted 2}\,, 
where since $v_5$ cannot reach $v_2$ with positive probability, 
it is removed from the system; this ensures that the Markov chain is irreducible.
Hence, it has a unique stationary distribution; this applies for all rerouted Markov chains of any finite absorbing Markov chain.
\end{example}

The following lemma will be useful for the upcoming result, \cref{theorem: hitting probability absorbing}.

\begin{lemma}\label{lemma: hitting probability absorbing helper deviation}
Let $a$ and $b$ be positive real numbers and let $a'$ and $b'$
be nonnegative real numbers.
For a given $\varepsilon$ satisfying $0 \leq \varepsilon < \frac{a+b}{2}$,
if $\left|a'-a\right|\leq \varepsilon$ and
$\left|b'-b\right| \leq \varepsilon$, then
\begin{align*}
    \frac{a-\varepsilon}{a+b} \leq \frac{a'}{a'+b'}
    \leq \frac{a+\varepsilon}{a+b}.
\end{align*}
\end{lemma}
\begin{proof}
Let $a'=a+\delta_1$ and $b'=b+\delta_2$,
where $|\delta_1|,|\delta_2| \leq \varepsilon$.
First, notice that $\frac{a'}{a'+b'}$ is well-defined,
as $a'+b' = a+b+\delta_1+\delta_2 \geq a+b-2\varepsilon > 0$.
Since $a' \geq 0$, $\frac{a'}{a'+b'}\geq 0$.
Holding $a$ and $b$ fixed, we see that the quantity 
$\frac{a'}{a'+b'} = 1/(1+\frac{b'}{a'})
= 1/(1+\frac{b+\delta_2}{a+\delta_1})$
is weakly increasing with respect to $\delta_1$
and weakly decreasing with respect to $\delta_2$.
(Technically this argument requires $a+\delta_1 > 0$
and therefore does not handle the case $a'=0$
but this case is easily dealt with separately.)
Hence, within the bounds for $\delta_1$ and $\delta_2$,
$\frac{a'}{a'+b'}$ is minimized
when $\delta_1 = -\varepsilon$ and $\delta_2 = \varepsilon$, yielding
\begin{align*}
    \frac{a'}{a'+b'}
    = \frac{a+\delta_1}{a+b+\delta_1+\delta_2}
    \geq \frac{a-\varepsilon}{a+b-\varepsilon+\varepsilon}
    = \frac{a-\varepsilon}{a+b},
\end{align*}
and $\frac{a'}{a'+b'}$ is maximized
when $\delta_1 = \varepsilon$ and $\delta_2 = -\varepsilon$, yielding
\begin{align*}
    \frac{a'}{a'+b'}
    = \frac{a+\delta_1}{a+b+\delta_1+\delta_2}
    \leq \frac{a+\varepsilon}{a+b+\varepsilon-\varepsilon}
    = \frac{a+\varepsilon}{a+b},
\end{align*}
which completes the proof.
\end{proof}

The following theorem demonstrates how 
the hitting probabilities of a finite absorbing Markov chain 
are approximated by the firing vector of the rerouted Markov chain.

\begin{theorem}\label{theorem: hitting probability absorbing}
Given a finite absorbing Markov chain, let $\v^{(N)}$ be the firing vector 
after $N$ steps of a hunger game process 
on the rerouted Markov chain associated with state $v$.
Then the sequence $\{a_N\}$ defined by
\begin{align*}
    a_N = \frac{\v^{(N)}_u}{\displaystyle\sum_{u' \in U} \v^{(N)}_{u'}},
\end{align*}
where we define $a_N$ to be 0 when the denominator equals 0, 
converges to the hitting probability $h_u(v)$ with discrepancy $O(1/N)$; 
that is, there exists a constant $C$ such that 
$a_N$ differs from $h_u(v)$ by at most $\frac{C}{N}$ for all $N$. 
\end{theorem}
\begin{proof}
In the rerouted Markov chain, 
$v$ can reach every state in $V$, and every state in $V$ 
can reach an absorbing state (as the original Markov chain was absorbing), 
which in the rerouted Markov chain moves back to $v$ with probability 1; 
as a result, the rerouted Markov chain is irreducible.
By \cref{theorem: irreducible fidelity converge stationary}\,,
the normalized firing vector $\frac{1}{N}\v^{(N)}$ 
converges to the unique stationary distribution $\ppi$ 
of the rerouted Markov chain within distance $\frac{c}{N}$ 
for some constant $c$.

It is a standard fact that the expected number of visits to state $w$ 
before returning to $v$ is given by $\frac{\ppi_w}{\ppi_v}$; 
see, for example, \cite[Theorem 1.7.6]{norris1998markov}.
But, as after visiting some $u'\in U$ the next state visited must be $v$,
meaning $u'$ can be visited at most once before returning to $v$,
the expected number of visits to any such $u'$ 
equals the probability of visiting $u'$ before returning to $v$.
As a result, for all $u \in U$, the hitting probability $h_{u}(v)$ 
is proportional to $\ppi_{u}$, and given that an absorbing Markov chain 
reaches an absorbing state with probability 1, we have
\begin{align*}
    h_{u}(v) = \frac{\ppi_{u}}{\displaystyle\sum_{u' \in U} \ppi_{u'}}.
\end{align*}

As $\sum_{u' \in U} \ppi_{u'}$ is a positive constant, 
there exists a finite $M$ such that for all $N>M$ we have
\begin{align*}
    \frac{c}{N} < \frac{1}{2} \sum_{u' \in U} \ppi_{u'}
\end{align*}
with $c$ as above.
As $a_N \in [0,1]$ is bounded, there exists a constant $C_1$ 
such that $a_N$ is within $\frac{C_1}{N}$ of $h_u(v)\in[0,1]$ for all $N \leq M$;
in particular, $C_1=M$ suffices.

For $N > M$, notice that
\begin{align*}
    a_N = \dfrac{\v_u^{(N)}}{\displaystyle\sum_{u'\in U} \v_{u'}^{(N)}}
    = \dfrac{\frac{1}{N}\v_u^{(N)}}{\frac{1}{N}\v_u^{(N)} + {\displaystyle\sum_{u'\in U\setminus\{u\}}} \frac{1}{N} \v_{u'}^{(N)}}.
\end{align*}
As $\frac{1}{N}\v^{(N)}$ is within distance $\frac{c}{N}$ of $\ppi$ in the $L^1$ metric, we know both $\left|\frac{1}{N}\v_u^{(N)}-\ppi_u\right| \leq \frac{c}{N}$ and
\begin{align*}
    \left|\displaystyle\sum_{u'\in U\setminus\{u\}}\frac{1}{N} \v_{u'}^{(N)} - \displaystyle\sum_{u' \in U\setminus\{u\}} \ppi_{u'}\right| \leq \frac{c}{N}.
\end{align*}
Due to the rerouted Markov chain being irreducible, we also know that $\ppi_v > 0$ for all states $v\in V$; in addition, each component of the visit vector $\v^{(N)}$ is a nonnegative integer.
Lastly, by definition of $M$, we have
\begin{align*}
    0 \leq \frac{c}{N} < \frac{1}{2} \left(\ppi_{u} + \sum_{u'\in U\setminus\{u\}} \ppi_{u'}\right).
\end{align*}
Letting $\varepsilon = \frac{c}{N}$, we apply
\cref{lemma: hitting probability absorbing helper deviation} with the values
\begin{align*}
    a = \ppi_u > 0, \qquad\quad
    b = \sum_{u'\in U\setminus\{u\}} \ppi_{u'} > 0, \qquad\quad
    a' = \frac{1}{N}\v_u^{(N)} \geq 0, \qquad\quad
    b' = \sum_{u'\in U\setminus\{u\}}\frac{1}{N} \v_{u'}^{(N)} \geq 0
\end{align*}
to find
\begin{align*}
    \frac{\ppi_{u} - \frac{c}{N}}{\ppi_{u} + \displaystyle\sum_{u'\in U\setminus\{u\}} \ppi_{u'}}
    \leq \dfrac{\frac{1}{N}\v_u^{(N)}}{\frac{1}{N}\v_u^{(N)} + {\displaystyle\sum_{u'\in U\setminus\{u\}}} \frac{1}{N} \v_{u'}^{(N)}}
    \leq \frac{\ppi_u + \frac{c}{N}}{\ppi_{u} + \displaystyle\sum_{u'\in U\setminus\{u\}} \ppi_{u'}},
\end{align*}
or equivalently
\begin{align*}
    h_u(v) - \frac{1}{N}\cdot\frac{c}{\displaystyle\sum_{u'\in U}\ppi_{u'}}
    \leq a_N
    \leq h_u(v) + \frac{1}{N}\cdot\frac{c}{\displaystyle\sum_{u'\in U}\ppi_{u'}}.
\end{align*}

Thus, the deviation of $a_N$ from $h_u(v)$ is at most
\begin{align*}
    \dfrac{1}{N}\cdot\frac{c}{\displaystyle\sum_{u' \in U} \ppi_{u'}}.
\end{align*}
The second factor is a constant, which we will denote $C_2$, 
and thus taking $C=\max(C_1,C_2)$ 
we find that for all $N$, $a_N$ differs from $h_u(v)$ by at most $\frac{C}{N}$.
\end{proof}

\begin{remark}\label{remark: escape probability}
We can apply \cref{theorem: hitting probability absorbing} 
to deterministically approximate the escape probability of $u$ from $v$ 
on an irreducible Markov chain, or the probability that a chain 
starting at state $v$ reaches state $u$ before returning to $v$.
We do this by splitting vertex $v$ into vertices $v_0$ and $v_1$ 
where all outgoing edges of $v$ now emanate from $v_0$ 
and all inbound edges of $v$ now end at $v_1$.
Furthermore, we may remove all outgoing edges from vertex $u$.
The resulting Markov chain is absorbing 
with the set of absorbing states $U=\{u,v\}$, 
and thus one may use \cref{theorem: hitting probability absorbing} 
on the rerouted Markov chain associated with state $v_0$ 
to approximate the hitting probability $h_{v_1}(v_0)$.
Notice that $h_{v_1}(v_0)$ equals the escape probability 
of $u$ from $v$ in the original irreducible Markov chain, 
and thus \cref{theorem: hitting probability absorbing} 
enables us to calculate escape probabilities.
\end{remark}

We now consider the \textbf{absorption time} of an absorbing Markov chain
whose set of absorbing states is $U \subseteq V$, given by
\begin{align*}
    T_U = \min\{t \geq 0: X_t \in U\}.
\end{align*}
At the absorption time, the chain enters an absorbing state, 
where it remains forever.

For example, starting at state $v_2$ in the absorbing Markov chain given in \cref{subfig: ex rerouted base}, with probability $\frac{1}{2}$ we enter absorbing state $v_1$ in the first step, and otherwise move to state $v_3$, which has the same 50-50 split between moving to an absorbing state ($v_4$) or not ($v_2$).
Thus $T_U$ in this case is a geometric random variable with success probability $\frac{1}{2}$, which has expected value 2.

The expected absorption time is $\E_v T_U$, 
where $\E_v$ denotes the expected value for the law 
of the Markov chain initialized at state $v$.
Then, when using the rerouted Markov chain, 
a chain can naturally be divided into epochs 
each of which ends with an occurrence of a state in $U$, 
which is then rerouted in the next time step to $v$.
Hence, the expected absorption time $\E_v T_U$ 
can be calculated from the stationary distribution $\ppi$ 
of the rerouted Markov chain associated with $v$ to be
\begin{align*}
    \E_v T_U = \frac{1}{\displaystyle\sum_{u \in U} \ppi_u}-1.
\end{align*}
Notice that if $v \in U$, then $\E_v T_U = 0=\frac{1}{1}-1$; 
the $-1$ term corresponds to the extra time step needed to reroute 
from an absorbing state to $v$ for each epoch.

The following theorem shows how the expected absorption time 
of an absorbing Markov chain can be approximated 
using the normalized firing vector of the rerouted Markov chain.

\begin{theorem}\label{theorem: absorption time}
Given an absorbing Markov chain, let $\v^{(N)}$ be the firing vector 
after $N$ steps of the hunger game process on the rerouted Markov chain 
associated with state $v$.
Then the sequence $\{b_N\}$ defined by
\begin{align*}
    b_N = \frac{N}{\displaystyle\sum_{u \in U} \v^{(N)}_u}-1,
\end{align*}
where we define $b_N$ to be 0 when the denominator equals 0, 
converges to the expected absorption time $\E_v T_U$, 
where there exists a constant $C$ such that $b_N$ 
is within distance $\frac{C}{N}$ of $\E_v T_U$ for all $N$.
\end{theorem}
\begin{proof}
From the proof of \cref{theorem: hitting probability absorbing}\,,
we see that the rerouted Markov chain 
has a unique stationary distribution $\ppi$, and thus 
the normalized firing vector $\frac{1}{N}\v^{(N)}$ 
converges to $\ppi$ within distance $\frac{c}{N}$ for some constant $c$.

As $\sum_{u \in U} \ppi_{u}$ is a positive constant, 
there exists a finite $M$ such that for all $N>M$ we have
\begin{align*}
    \frac{c}{N} < \frac{1}{2} \sum_{u \in U} \ppi_{u}.
\end{align*}
As $b_N$ is finite, there exists a constant $C_1$ 
such that $b_N$ is within $\frac{C_1}{N}$ of $\E_v T_U$ for all $N \leq M$.
For $N > M$, it suffices to show the existence of a constant $C_2$ 
such that for all $N > M$,
\begin{align*}
    \left| \frac{1}{\displaystyle\sum_{u \in U} \frac{1}{N}\v^{(N)}_u} - \frac{1}{\displaystyle\sum_{u \in U} \ppi_u} \right| < \frac{C_2}{N}. 
\end{align*}
The furthest away $b_N$ can be from $\E_v T_U$ 
results in the left hand side becoming
\begin{align*}
    \frac{1}{-\displaystyle\frac{c}{N} + \displaystyle\sum_{u \in U} \ppi_u} - \frac{1}{\displaystyle\sum_{u \in U} \ppi_u} = \frac{1}{N}\cdot\frac{c}{\left(-\displaystyle\frac{c}{N} + \displaystyle\sum_{u \in U} \ppi_u\right)\left(\displaystyle\sum_{u \in U} \ppi_u\right)}.
\end{align*}
As $N>M$, we have
\begin{align*}
    -\frac{c}{N}+\displaystyle\sum_{u \in U} \ppi_{u} > \frac{1}{2}\sum_{u \in U} \ppi_{u} > 0,
\end{align*}
so in the worst case the deviation from $\E_v T_U$ is less than
\begin{align*}
    \frac{1}{N}\cdot \frac{2c}{\left(\displaystyle\sum_{u \in U} \ppi_{u}\right)^2}.
\end{align*}
The second factor is a constant, which we will denote $C_2$, 
and thus taking $C=\max(C_1,C_2)$ yields that $b_N$ 
is at most $\frac{C}{N}$ away from $\E_v T_U$ for all $N$.
\end{proof}

\begin{remark}\label{remark: expected return time irreducible}
We can apply \cref{theorem: absorption time} to calculate 
the expected return time for state $v$ in an irreducible Markov chain.
We do this in a similar method to \cref{remark: escape probability} 
by splitting vertex $v$ into vertices $v_0$ and $v_1$ 
where all outgoing edges of $v$ now emanate from $v_0$ 
and all inbound edges of $v$ now end at $v_1$.
The resulting Markov chain is absorbing with unique absorbing state $v_1$, 
and thus one may use \cref{theorem: absorption time} 
on the rerouted Markov chain associated with state $v_0$ 
to approximate the expected absorption time.
As $v_1$ is the only absorbing state, this calculates 
the expected hitting time from $v_0$ to $v_1$, 
which is the expected return time from $v$ 
in the original irreducible Markov chain.
\end{remark}

%% file: recurrence.tex
\section{Recurrent states and the basin of attraction}\label{section: recurrence}

In this section we restrict our attention to irreducible finite Markov chains
with no absorbing state.
Additionally we assume that all transition probabilities are rational.

From the argument in \cref{remark: finite orbit}, 
all hunger vectors $\h\in\R^n$ are pre-periodic (i.e.,
eventually enter a cycle).
We say that a vector that is part a cycle
(i.e., returns to itself after a finite number of steps) is \textbf{recurrent},
and we define the \textbf{basin of attraction} to be 
the set of recurrent vectors on the hyperplane 
$Z=\{\h\mid \h\cdot \mathbf{1} = 0\}$ 
(the hyperplane of hunger vectors with total hunger 0).
The basin of attraction, being the set of recurrent vectors, 
is thus in some sense analogous to the critical group $\Z^n/\Z^n\Delta$ 
in the context of chip-firing.

\begin{example}\label{example: basin of attraction 1}
Consider the finite rational irreducible Markov chain with hunger matrix
\[ H = \begin{bmatrix} -0.2 & 0.2 & 0 \\ 0.2 & -0.6 & 0.4 \\ 0.6 & 0.4 & -1 \end{bmatrix},\]
whose unique stationary distribution is $\left(\frac{11}{18},\frac{5}{18},\frac{2}{18}\right)$.
Its basin of attraction is given in \cref{fig: basin of attraction 1}.
The basin is partitioned into four colors based on the firing order, up to cyclic shift.
Red corresponds to the cyclic firing order 231121121311211211, i.e., firing state 2, then state 3, then state 1, and so on.
Green corresponds to 312112112311211211, cyan corresponds to 231112112311211211, and violet corresponds to 312112112311121121.
As the visit vector for each of these 18-step cycles must satisfy $\v H = \mathbf{0}$, i.e., be a multiple of the stationary distribution, we find each cyclic firing order consists of state 1 firing eleven times, state 2 firing five times, and state 3 firing two times.
\begin{figure}[htb]
    \centering
    \includegraphics[width=0.5\textwidth]{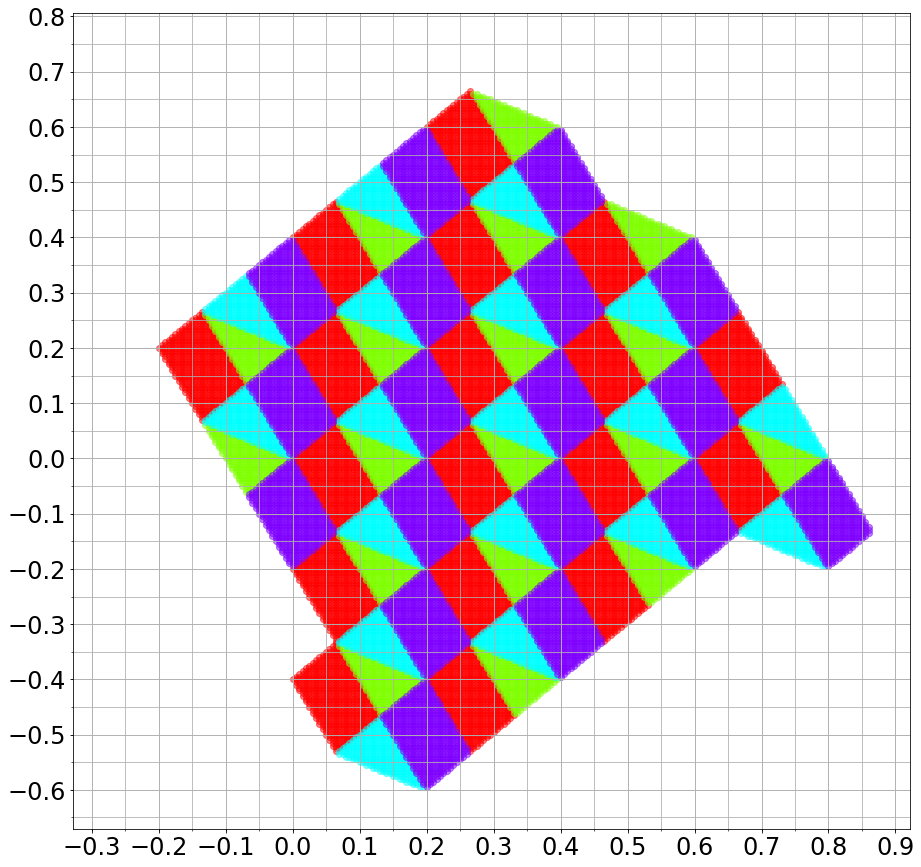}
    \caption{The basin of attraction for a finite rational irreducible Markov chain with three states.
    The basin lies in the hyperplane $x + y + z = 0$, so the third coordinate is omitted in order to plot it on the $x,y$-plane.
    For example, the point $(0.8,-0.1)$ in the rightmost violet region corresponds to $\h = (0.8,-0.1,-0.7)$.
    The basin is partitioned into four colors based on the order of firings; for example, the vectors in the violet portions of the basin all follow the same sequence of firings, up to cyclic shift, i.e., choosing a starting index for the periodic firing sequence.
    So, the four colors correspond to four distinct cyclic firing orders.}
    \label{fig: basin of attraction 1}
\end{figure}
\end{example}

\begin{remark}\label{remark: convex partition basin}
Numerous observations can be made from \cref{example: basin of attraction 1} 
and \cref{fig: basin of attraction 1}.
Firstly, the partitioning based on cyclic firing order appears to yield 
discrete congruent pieces; in \cref{fig: basin of attraction 1}, 
the red and violet colors are split into congruent parallelograms, 
and the cyan and green colors are each split into congruent triangles.
Each color corresponds to a cyclic firing order, and if this cycle has period $p$, 
then the set of vectors of this particular color 
can be partitioned into $p$ sets based on the firing order of their cycle, 
no longer up to cyclic shift.
Each of these sub-pieces are mapped to each other in a cycle 
under the hunger game process, 
so clearly the sub-pieces of a given color are all congruent.

While evidently the entire basin of attraction need not be convex, 
each sub-piece of the basin of attraction, e.g., 
the parallelograms and triangles of \cref{fig: basin of attraction 1}, is convex; 
this holds for all finite rational Markov chains.
To see why this is the case, notice that each sub-piece 
corresponds to a unique firing order that returns a hunger vector to itself.
The set of vectors that will fire to a given state next under the hunger game process 
is given by the intersection of various linear inequalities, 
namely the inequalities that ensure this state had the highest hunger.
Hence, the set of vectors that follow the firing order corresponding to a given sub-piece 
is given by the intersection of various linear inequalities, 
one set of inequalities for each firing step of the hunger game process.
As the intersection of linear inequalities is convex, we find each sub-piece is convex.
\end{remark}

Intuitively, the basin of attraction should be near the origin, 
for the hunger game naturally attempts to equilibrate the hunger vector 
so that each state has approximately the same hunger.
The following result supports this intuition by demonstrating that the origin 
is always in the basin of attraction for any finite rational irreducible Markov chain.

\begin{proposition}\label{proposition: 0 hunger recurrent}
For any finite rational irreducible Markov chain, 
the zero vector $\mathbf{0}$ is a recurrent hunger state.
Moreover, the number of steps needed to return to $\mathbf{0}$ 
is the least common denominator of the stationary probabilities of the Markov chain.
\end{proposition}
\begin{proof}
Let $d$ be the least common denominator 
of the transition probabilities of the Markov chain.
Viewing the Markov chain as a weighted digraph $G$
with each vertex having outgoing weights summing to 1,
we can multiply all weights by $d$ so that all weights are now positive integers, 
and convert weighted edges into multiedges to obtain a directed multigraph $G'$
in which every vertex has outdegree $d$.

A finite irreducible Markov chain with rational transition probabilities
has a unique stationary distribution $\ppi$, all of whose entries are rational.
Let $n$ be the least common denominator of $\ppi$, so that $n\ppi$
is the unique multiple of $\ppi$ 
that has nonnegative, mutually coprime integer entries, and let $\p = n\ppi$.
Since the sum of the entries of $\ppi$ is 1, the sum of the entries of $\p$ is $n$. 
As $\ppi H=\mathbf{0}$, we have $\p H=\mathbf{0}$.

Starting at $\mathbf{0}$, consider the first vertex $v$ that fires more than $\p_v$ times, 
and consider the hunger state just before this vertex fires for the $(\p_v+1)$-th time.
Let the firing vector at this step be $\x$, so that $\x_u$ is the number of times 
state $u$ received the chip, and thus $\x H$ is the current hunger state.
Because $\x_v = \p_v$ and 
because for all $u \neq v$ state $u$ has fired at most $\p_u$ times, 
$\w := \p - \x$ satisfies $\w_v = 0$ and $\w_u \geq 0$ for all $u \neq v$.
Firing vertices $u\neq v$ can only (weakly) increase the hunger of $v$, 
so $(\w H)_v \geq 0$.
To prove that last assertion more formally, recall that $\w_v = 0$, 
so we can write $\w=\sum_{u \neq v} c_u \e^{(u)}$ with $c_u \geq 0$ for all $u$,
where $\e^{(u)}$ is the elementary basis vector with a 1 at index $u$ and 0 otherwise.
Then we have
\begin{align*}
    (\w H)_v
    &= \sum_{u \neq v} c_u (\e^{(u)} H)_v
    = \sum_{u \neq v} c_u H_{uv}
    = \sum_{u \neq v} c_u P_{uv}
    \geq 0,
\end{align*}
as claimed. This yields that 
\begin{align*}
    (\x H)_v
    = (\p H - \w H)_v
    = -(\w H)_v
    \leq 0.
\end{align*}
Yet with a hunger state of $\x H$, vertex $v$ received the chip, 
meaning it has the highest hunger.
This means every state has hunger at most $(\x H)_v \leq 0$, so total hunger is at most 0.
However, total hunger is invariant under the hunger game process;
since the initial hunger vector was $\mathbf{0}$ with total hunger 0, 
total hunger must still be equal to 0.
This equality case requires that $\x H = \mathbf{0}$, 
implying $\x$ is a multiple of $\ppi$. 
Since $\x_v = \p_v$ and $\p = n\ppi$, we deduce $\x = \p$.

Thus the hunger state is $\x H=\mathbf{0}$, so $\mathbf{0}$ is recurrent.
Moreover, the number of steps needed to return to $\mathbf{0}$ 
is the sum of the entries of $\x=\p=n\ppi$, which is $n$, as claimed.
\end{proof}

From empirical observations, we conjecture that 
the periods of all cycles in the hunger game for a given chain are equal, 
and specifically are equal to the $n$ 
that was shown in \cref{proposition: 0 hunger recurrent}
to be the period of the hunger vector $\mathbf{0}$.

\begin{conjecture}\label{conjecture: period stationary lcd}
The period of every cycle in the hunger game 
for a given finite rational irreducible Markov chain 
is the least common denominator of the stationary probabilities of the Markov chain.
\end{conjecture}

Visually, \cref{conjecture: period stationary lcd} 
applied to the Markov chain from \cref{example: basin of attraction 1} 
corresponds to the observation that each color in \cref{fig: basin of attraction 1} 
consists of the same number of congruent pieces, namely 18.

\begin{figure}[htbp]
    \centering
    \includegraphics[width=0.6\textwidth]{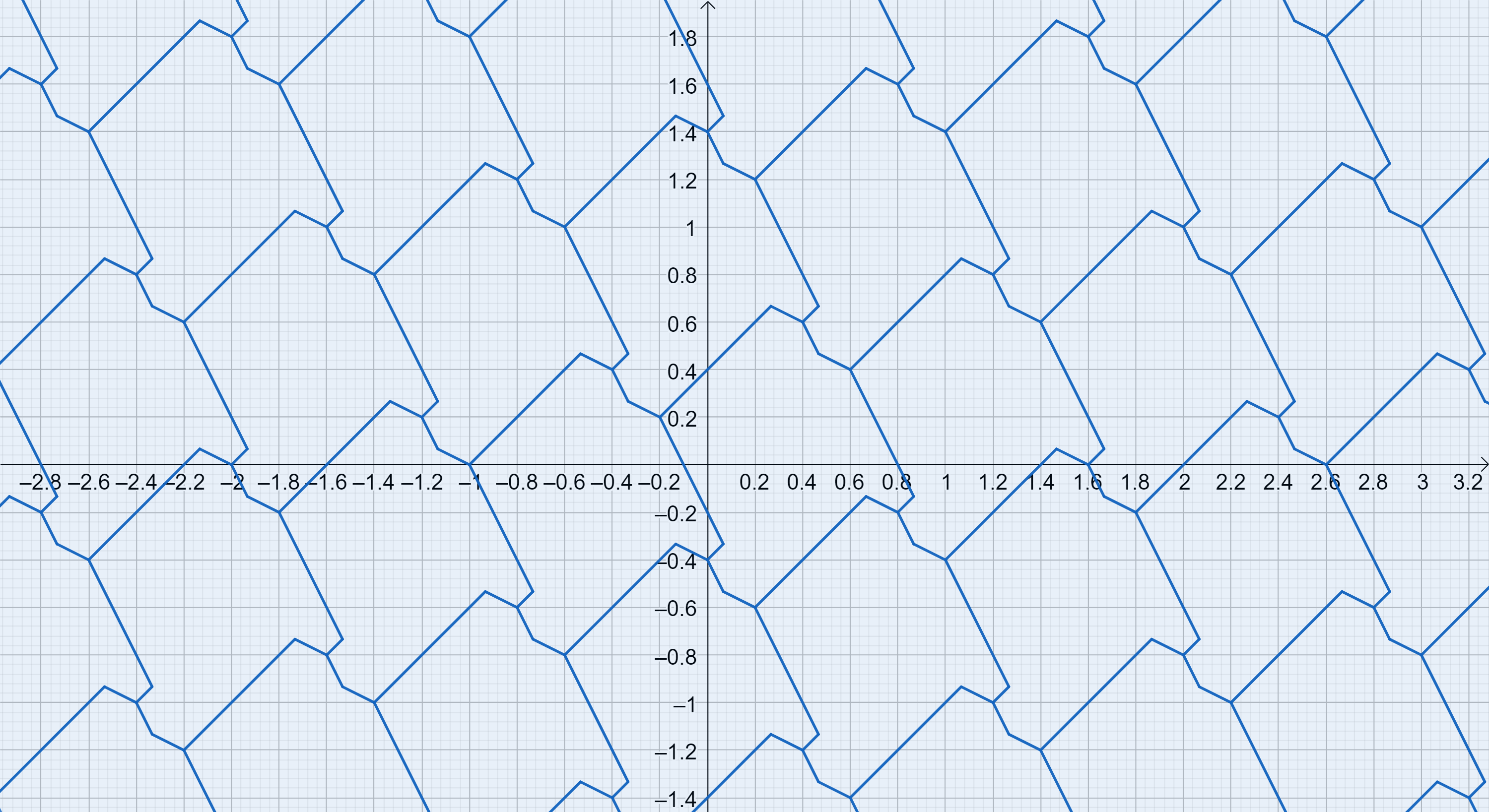}
    \caption{The basin of attraction from \cref{example: basin of attraction 1} tiles the (hyper)plane $Z$.}
    \label{fig: basin tiling}
\end{figure}

The basin of attraction in \cref{fig: basin of attraction 1} tiles the hyperplane $Z$, as illustrated in \cref{fig: basin tiling}; 
in particular, one can observe that its concave boundaries fit complementarily 
with the opposite side of the basin of attraction.
The following conjecture formalizes this observation 
and poses it for general Markov chains.

\begin{conjecture}\label{conjecture: basin tiles hyperplane}
For a finite rational irreducible Markov chain, 
the basin of attraction tiles the hyperplane $Z$ by translation, 
and the translation vectors that relate each tile with each other 
forms an $(n-1)$-dimensional sublattice in $Z$ of $\R^n$ 
with basis vectors $H_1-H_2,H_2-H_3,\dots,H_{n-1}-H_n$.
\end{conjecture}

In fact, \cref{conjecture: basin tiles hyperplane} 
implies \cref{conjecture: period stationary lcd}, as seen in the following proposition.

\begin{proposition}\label{proposition: tiling conjecture implies periodicity conjecture}
For a finite rational irreducible Markov chain, 
if the basin of attraction tiles the hyperplane $Z$ 
with lattice of translation vectors generated by $H_1-H_2,\dots,H_{n-1}-H_n$, 
then the period of any cycle under the hunger game process is 
the least common denominator of the stationary probabilities of the Markov chain.
\end{proposition}
\begin{proof}
By using \cref{proposition: 0 hunger recurrent}, it suffices to show that 
the tiling condition implies that all periods in the basin of attraction are the same.
Suppose the lattice in the statement of the proposition is $L$.
We say that two points in $Z$ are equivalent mod $L$ if they differ by 
an element of $L$, and we write the set of equivalence classes as $Z/L$.
Now view the hunger game process, which naturally acts upon $Z$, 
as acting instead on $Z/L$, where the tiling condition implies that 
the basin of attraction can serve as a set of coset representatives for $Z/L$.
Initially this might seem like nonsense since the choice of which state to fire
depends on inequalities relating the elements of the hunger vector 
and these inequalities get washed out when we mod out by $L$,
but all rows of $H$ are equivalent mod $L$, 
i.e., each row of $H$ has the same image under the projection map $Z\to Z/L$,
so the choice of which state gets fired is moot;
we may as well suppose that the $H_1$ coset
is added at each stage, as if state 1 were firing repeatedly.
Using the basin of attraction as a set of coset representatives 
and observing that a vector inside the basin 
must stay within the basin under the hunger game process, 
we find the hunger game process on $Z/L$ is equivalent to 
the original hunger game on $Z$ for any vector in the basin.
If the zero vector $\0$ has period $p$, then $pH_1=0$ in $Z/L$.
Then let $q$ be the minimum of the set of positive integers $q'$ 
such that $q'H_1=0$ in $Z/L$.
After $q$ steps of the hunger game process, 
every vector in the basin thus returns to itself, 
so we find $p=q$ and every vector in the basin of attraction must have period $p$, 
which completes the proof.
\end{proof}

We provide a partial result towards \cref{conjecture: basin tiles hyperplane}, 
proving that the basin of attraction translated under the stated lattice 
covers the hyperplane $Z$; 
to prove \cref{conjecture: basin tiles hyperplane} 
and thus \cref{conjecture: period stationary lcd}, 
it remains to show this covering has no overlap.

\begin{proposition}\label{proposition: basin covers hyperplane}
For a finite rational irreducible Markov chain, 
let $L\subset Z$ be the sublattice in $Z$ of $\R^n$ 
with basis vectors $H_1-H_2,\dots,H_{n-1}-H_n$.
Then for any vector $\h\in Z$, there exists a $\u\in L$ 
such that $\h-\u$ lies in the basin of attraction.
\end{proposition}
\begin{proof}
As $\h$ is pre-periodic, after some finite number of steps 
it will reach some vector in the basin, say after $t_0$ steps.
Additionally, every vector in the basin 
stays in the basin under the hunger game process, 
so for all integers $t \geq t_0$, 
applying the hunger game process $t$ times to $\h$ yields a vector in the basin.
Let $p$ be the least common denominator of the stationary probabilities 
of the unique stationary distribution $\ppi$ of our Markov chain, 
and fix $t$ to be the minimum nonnegative integer $\geq t_0$ such that 
$t$ is a multiple of $p$, say $kp$ for positive integer $k$.
Then applying the hunger game process $t$ times to $\h$ 
yields a vector $\x \in Z$ in the basin of attraction, 
say after firing state $i$ a total of $\v_i$ times for each $i$.
Constructing vector $\v$ from these values, 
as $\x$ was reached after firing $\h$ a total of $t$ times, 
we have $\x=\h+\v H$ and $\v\cdot\1=t=kp$.
Notice that $p\ppi$ is the primitive integer vector 
in the direction of the stationary distribution, so we have $p\ppi H= \0$.
Hence $\x = \h + (\v - kp\ppi)H$, where $(\v-kp\ppi)\cdot \1 = 0$.
As $L=\{\w H\mid\w\cdot\1=0\}$, 
letting $\u=(kp\ppi-\v)H\in L$ yields $\x = \h - \u$, as desired.
\end{proof}

%% file: conclusion.tex
\section{Comments}\label{section: conclusion}
The hunger game achieves high fidelity (that is, low discrepancy)
for the frequency with which a specified state occurs,
but not for the frequency with which two specified states occur in succession;
the rotor-router game achieves fidelity for the frequency with which
two specified states occur in succession, but not for the frequency 
with which three specified states occur in succession.
One could use a block-encoding trick
(with new states encoding pairs of old states) 
to construct a simulation scheme that
ensures fidelity for the frequency with which
three specified (old) states occur in succession, but at a price:
the resulting rotor-router network will have more rotors,
and hence the quantities that upper-bound discrepancy
in \cite{holroyd2010rotor} will be larger.
Indeed, consider derandomization of a fair coin process.
De Bruijn \cite{debruijn1946} showed that 
for each $k$, there is a cyclic word of length $2^k$
that contains each possible bit-string of length $k$ exactly once;
it follows that if we let $W$ be 
the infinite word associated with this cyclic word,
then for each bit-string $w$ of length $k$,
the unnormalized discrepancy $D(w,n)$,
defined as $n/2^k$ minus the number of apperances of $w$
in the first $n+k-1$ positions of $W$, is bounded.
On the other hand, one cannot find a single infinite word $W$
that accomplishes this task simultaneously for all $k$,
since for example if $w$ is a string containing one or more 1s,
the boundedness of $D(w,n)$
is inconsistent with the occurrence of arbitrarily long strings of 0's.

We propose that a sort of trade-off principle is at work here: 
the more numerical characteristics of a random process
that we attempt to mimic in a deterministic simulation,
the worse our mimicry will be, if our criterion of merit is the rate at which
the simulation's numerical characteristics converge to the limiting values 
of the corresponding numerical characteristics of the random process.
In this context, one sometimes sees the words
``psuedorandom'' and ``quasirandom'',
where a pseudorandom simulation is expected to pass ``all'' tests of randomness
while a quasirandom simulation is expected only to pass certain specified tests.
It is reasonable to expect that the more asymptotic tests of randomness
a quasirandom simulation is required to pass,
the worse the convergence will be.
Indeed, for certain tests of randomness,
the previous assertion is almost a tautology;
for if (say) a random sequence of $0$s and $1$s
had its running average approach 1/2 with discrepancy smaller than $\pm \sqrt{N}$,
that rapidity of convergence would ipso facto
demonstrate negative correlation between the bits,
constituting failure of a correlational test of randomness.

Applying this trade-off principle in the reverse direction,
it seems reasonable to hope that by relaxing our stringency,
and requiring a quasirandom simulation scheme
to pass {\em fewer} tests for randomness,
we enable it to do a {\em better} job
with the tests that it is required to pass;
that is, we make faster convergence possible.
Hence, the delocalized firing mechanism of the hunger game,
by no longer preserving the frequency with which
two specified states occur in succession, may allow
the hunger game to outperform the rotor-router in convergence.

\begin{figure}[htbp]
    \centering
    \begin{tikzpicture}[scale=0.6,font=\normalsize,baseline,thick]
        \foreach \x in {1,...,2} {
            \filldraw[color=black,fill=green!30,thick] (4*\x,0) circle (16pt);
            \node at (4*\x,0) {$v_{\x}$};
        }
        \foreach \x in {2} {
            \draw [->,>=latex] plot [smooth,tension=1] coordinates {(4*\x-0.866*0.7,0.4*0.7) (4*\x-2,0.7) (4*\x-4+0.866*0.7,0.4*0.7)};
        }
        \foreach \x in {1} {
            \draw [->,>=latex] plot [smooth,tension=1] coordinates {(4*\x+0.866*0.7,-0.4*0.7) (4*\x+2,-0.7) (4*\x+4-0.866*0.7,-0.4*0.7)};
        }
        \node at (6,1) {$q$};
        \node at (6,-1) {$q$};
        \draw [->,>=latex] plot [smooth,tension=5] coordinates {(4-0.866*0.7,0.5*0.7) (4-1.5,0) (4-0.866*0.7,-0.5*0.7)};
        \node at (1.5,0) {$1-q$};
        \draw [->,>=latex] plot [smooth,tension=5] coordinates {(8+0.866*0.7,-0.5*0.7) (8+1.5,0) (8+0.866*0.7,0.5*0.7)};
        \node at (10.5,0) {$1-q$};
    \end{tikzpicture}
    \caption{A Markov chain where the hunger game outperforms the rotor-router in convergence}
    \label{fig: conclusion outperform}
\end{figure}
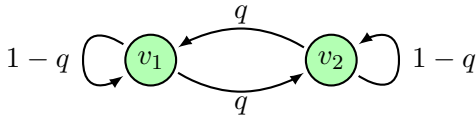

For example, consider the Markov chain with transition matrix
\begin{align*}
    \begin{bmatrix}
        1-q & q \\
        q & 1-q 
    \end{bmatrix}
\end{align*}
with $q$ small, as shown in \cref{fig: conclusion outperform}\,.
Then under rotor-routing the itinerary of the chip will be
periodic with repeating pattern $1,1,\dots,1,2,2,\dots,2$,
while under the hunger game the itinerary will simply alternate between 1 and 2,
yielding lower discrepancy for the visit-frequencies of the two states.